\documentclass[12pt]{article}

\usepackage{graphics,amsthm, amsmath, amssymb}

\usepackage{psfig}
\usepackage{color}

\expandafter\chardef\csname pre amssym.def at\endcsname=\the\catcode`\@
\catcode`\@=11

\def\undefine#1{\let#1\undefined}
\def\newsymbol#1#2#3#4#5{\let\next@\relax
 \ifnum#2=\@ne\let\next@\msafam@\else
 \ifnum#2=\tw@\let\next@\msbfam@\fi\fi
 \mathchardef#1="#3\next@#4#5}
\def\mathhexbox@#1#2#3{\relax
 \ifmmode\mathpalette{}{\m@th\mathchar"#1#2#3}%
 \else\leavevmode\hbox{$\m@th\mathchar"#1#2#3$}\fi}
\def\hexnumber@#1{\ifcase#1 0\or 1\or 2\or 3\or 4\or 5\or 6\or 7\or 8\or
 9\or A\or B\or C\or D\or E\or F\fi}

\font\tenmsa=msam10
\font\sevenmsa=msam7
\font\fivemsa=msam5
\newfam\msafam
\textfont\msafam=\tenmsa
\scriptfont\msafam=\sevenmsa
\scriptscriptfont\msafam=\fivemsa
\edef\msafam@{\hexnumber@\msafam}
\mathchardef\dabar@"0\msafam@39
\def\dashrightarrow{\mathrel{\dabar@\dabar@\mathchar"0\msafam@4B}}
\def\dashleftarrow{\mathrel{\mathchar"0\msafam@4C\dabar@\dabar@}}

\font\tenmsb=msbm10
\font\sevenmsb=msbm7
\font\fivemsb=msbm5
\newfam\msbfam
\textfont\msbfam=\tenmsb
\scriptfont\msbfam=\sevenmsb
\scriptscriptfont\msbfam=\fivemsb
\edef\msbfam@{\hexnumber@\msbfam}
\def\Bbb#1{\fam\msbfam\relax#1}

\theoremstyle{plain}
\newtheorem{theorem}{Theorem}[section]
\newtheorem{lemma}[theorem]{Lemma}

\theoremstyle{definition}
\newtheorem{definition}[theorem]{Definition}
\newtheorem{example}[theorem]{Example}

\newtheorem{remark}[theorem]{Remark}

\def\resultant{{\rm Res}}
\def\deg{{\rm deg}}
\def\degree{{\rm deg}}
\def\mul{{\rm mult}}
\def\emul{{\rm mult}_\epsilon}
\def\cP{{\mathcal P}}
\def\cc{{\mathcal C}}

\def\we{{\rm weight}}
\def\inn{{\rm in}}
\def\out{{\rm out}}

\def\sing{{\rm Sing}_{\epsilon}}
\def\cluster{{\mathfrak Cluster}}
\def\rad{{\rm radius}}

\def\qed{\hfill  \framebox(5,5){}}

\def\parasmall{\vspace{2 mm}}

\def\cPP{\overline{\cP}^{\Delta}(\Delta,t)}
\def\cSSy{{S}_{2}^{\Delta}(\Delta,y,t)}
\def\cBBy{B_{2}^{\Delta}(\Delta,y,t)}
\def\cRRy{R_{2}^{\Delta}(\Delta,y,t)}

\def\cPP{\overline{\cP}^{\Delta}(\Delta,t)}
\def\cSSx{{S}_{1}^{\Delta}(\Delta,x,t)}
\def\cBBx{B_{1}^{\Delta}(\Delta,x,t)}
\def\cRRx{R_{1}^{\Delta}(\Delta,x,t)}

\def\cPPsimple{\overline{\cP}^{\Delta}}
\def\cSSsimplex{{S}_{1}^{\Delta}}
\def\cBBsimplex{B_{1}^{\Delta}}
\def\cRRsimplex{R_{1}^{\Delta}}

\def\cPPsimple{\overline{\cP}^{\Delta}}
\def\cSSsimpley{{S}_{2}^{\Delta}}
\def\cBBsimpley{B_{2}^{\Delta}}
\def\cRRsimpley{R_{2}^{\Delta}}

\def\card{{\rm Card}}

\begin{document}

\title{Approximate Parametrization of Plane Algebraic Curves by Linear Systems of Curves\thanks{
Authors partially supported by   the Spanish "Ministerio de Educaci\'on y Ciencia" under the
Project MTM2005-08690-C02-01.}}


\author{Sonia P\'erez-D\'{\i}az and    J. Rafael  Sendra\\
Dpto. de Matem\'aticas \\
        Universidad de Alcal\'a \\
      E-28871 Madrid, Spain  \\
sonia.perez@uah.es, rafael.sendra@uah.es
\and Sonia L. Rueda  \\
Dpto. de Matem\'atica Aplicada \\
      E.T.S. Arquitectura,
 Universidad Polit\'ecnica de Madrid \\
        E-28040 Madrid, Spain \\sonialuisa.rueda@upm.es
\and Juana  Sendra \\
Dpto. de Matem\'atica Aplicada a la I.T. de Telecomunicaci\'on  \\
       E.U.I.T.Telecomunicaci\'on, Universidad Polit\'ecnica de Madrid \\
        E-28031 Madrid, Spain \\jsendra@euitt.upm.es
}
\date{}          
\maketitle

\begin{abstract}
It is well known that an irreducible algebraic curve is rational
(i.e. parametric) if and only if its genus is zero.   In this
paper, given a tolerance $\epsilon>0$ and an
$\epsilon$-irreducible algebraic affine plane curve $\cc$ of
proper degree $d$, we introduce the notion of
$\epsilon$-rationality, and we provide an algorithm to parametrize
approximately  affine $\epsilon$-rational plane curves, without
exact singularities at infinity, by means of linear systems of
$(d-2)$-degree curves. The algorithm outputs  a rational
parametrization of a  rational curve $\overline{\cc}$ of degree at
most $d$ which has the same points at infinity as $\cc$. Moreover,
although we do not provide a theoretical analysis, our empirical
analysis shows that $\overline{\cc}$ and $\cc$ are close in
practice.
\end{abstract}

\section*{Introduction}
Let $\mathfrak{O}^*$ be an algebraic or geometric object that
satisfies a property $\mathfrak{P}$ that implies the existence of
certain associated objects $\mathfrak{O}_{i}^*$; for instance,
$\mathfrak{O}^*$  might be a polynomial, $\mathfrak{P}$ the fact
of being reducible and $\mathfrak{O}_{i}^*$ the irreducible
factors. Computer algebra techniques provide, for a wide class of
situations, algorithms to check $\mathfrak{P}$, and to compute
exactly the associated objects $\mathfrak{O}_{i}^*$. However, in
many practical applications, we receive a perturbation
${\mathfrak{O}}$  of $\mathfrak{O}^*$, where $\mathfrak{P}$ does
not hold anymore neither the associated objects
$\mathfrak{O}_{i}^*$ exist. The problem, then, consists in
computing a new object $\overline{\mathfrak{O}}$, close to
$\mathfrak{O}$ and satisfying $\mathfrak{P}$, as well as the
associated objects $\overline{\mathfrak{O}}_i$ to
$\overline{\mathfrak{O}}$. We call approximate to an algorithm
solving a problem of the above type. Here, the notion of
``closeness" depends in general on the particular problem that one
is solving.

\parasmall

One can find in the literature approximate algorithms for computing gcds (see
 \cite{cor:gi:tra:wat}, \cite{em:ga:lo}, \cite{Pan1}), factoring polynomials (see \cite{Cor2}, \cite{GaRu}, \cite{pan},
\cite{Sa}), etc. For algebraic varieties there also exist
approximate solutions: see \cite{Cor4}, \cite{Dokken} for the
implicitization problem, in \cite{far:ra} the numerical condition
of implicitly given algebraic curves and surfaces has been
analyzed, and see \cite{BaRo}, \cite{Josef02}, \cite{Hart},
\cite{PSS}, \cite{PSS2} where the parametrization questions are
treated.

\parasmall

In this paper we consider the approximate parametrization problem
for affine plane algebraic curves. That is, with the above
terminology, $\mathfrak{O}^*$ is an  affine plane curve,
$\mathfrak{P}$ is the fact of being rational, and
$\mathfrak{O}^{*}_i$ is a rational parametrization of
$\mathfrak{O}^*$. So, the problem is stated as follows: we are
given an affine curve (say that it is a perturbation of a rational
curve) and we want to compute a rational parametrization of a
rational affine curve near  it; where we use the notion of
``vecinity" introduced in \cite{PSS}.

\parasmall

In \cite{PSS} and \cite{PSS2} the approximate parametrization
problem is solved for the special case of affine plane curves and
affine surfaces being a perturbation of a monomial curve and
surface, respectively. In both papers, the  basic tool is the use
of $\epsilon$-points (see also \cite{PSS3}). More precisely, given
a tolerance $\epsilon>0$, in \cite{PSS},  the parametrization
problem is solved for the case of affine plane curves having an
$\epsilon$-singularity of maximum multiplicity, and in \cite{PSS2}
the problem is solved for affine surfaces having also an
$\epsilon$-singularity of maximum multiplicity. The basic idea was
to use a pencil of lines through the $\epsilon$-singularity and,
hence, it was solved working as in the exact case for monomial
varieties.

\parasmall

In this paper, we generalize the ideas in \cite{PSS} to the case
of affine plane curves without singularities at infinity.  For
this purpose, the first obstacle is to associate suitably the
different $\epsilon$-singularities. This leads to the notion of
cluster. Then, we introduce the notion of (affine)
$\epsilon$-rationality, and we provide an algorithm to parametrize
approximately $\epsilon$-rational curves without exact
singularities at infinity. The idea of the algorithm is to work
with linear systems of curves of degree $d-2$, where $d$ is the
degree of the input curve. This system plays the role of the
linear system of adjoint curves in the exact parametrization
algorithm. In addition, we prove that the degree of the output
rational curve is bounded by the degree of the input one, and that
both curves have the same points at infinity. Differently to
\cite{PSS} we do not provide a theoretical analysis of the error
(i.e. on the closeness of input and output). However, our
empirical analysis shows that the curves are in practice near, and
it allows us to think about a theoretical treatment of this fact
as a future project.

\parasmall

The paper is structured as follows. In Section 1 we recall the
main notions and properties on $\epsilon$-singularities. Section 2
is devoted to recall the main ideas of the exact parametrization
algorithm for curves. In Section 3  we develop the idea of cluster
and we introduce the notion of  $\epsilon$-rationality. In Section
4 we derive the approximate algorithm, as well as the main
properties of the output curve. In Section 5 we illustrate the
algorithm by some example, and in Section 6 we analyze empirically
the error.

\parasmall

Throughout this paper, we use the following {\bf terminology}.
$\|\cdot\|$  and $\|\cdot \|_2$ denote the polynomial  $\infty$--norm and the usual unitary norm in ${\Bbb C}^2$, respectively. $|\cdot|$ denotes the module in $\Bbb C$.
The partial derivatives of a polynomial $g\in {\Bbb C}[x,
y]$  are denoted by
$g^{\overrightarrow{v}}:=\frac{\partial^{i+j} g}{\partial^{i} x
\partial^{j} y}$ where $\overrightarrow{v}=(i,j)\in {\Bbb N}^{2}$;
we assume that $g^{\tiny{\overrightarrow{0}}}=g.$ Moreover, for
$\overrightarrow{v}=(i,j)\in {\Bbb N}^{2}$,
$|\overrightarrow{v}|=i+j$. Also,   $\overrightarrow{e_1}=(1,0)$, and $\overrightarrow{e_2}=(0,1)$.

\parasmall

 In addition, we use the following  {\bf general assumptions}.
A tolerance $\epsilon$  is fixed such that $0 <\epsilon <1$. $\cal
C$ is an affine real plane algebraic curve over $\Bbb C$ of proper
degree $d>0$ (see Def. \ref{def-proper-degree}), without (exact)
singularities at infinity, not passing through $(1:0:0), (0:1:0)$,
and     defined by an $\epsilon$-irreducible polynomial $f(x,y)\in
{\Bbb R}[x,y]$ ; that is $f$ can not be expressed as $f(x,
y)=g(x,y)h(x,y)+{\cal E}(x, y)$ where $h,g,{\cal E}\in {\Bbb C}[x,
y]$ and $\|{\cal E}(x, y)\|<\epsilon \|f(x, y)\|$ (see
\cite{Cor2}, \cite{Kaltofen}). We denote by $\cc^h$ the projective closure of
$\cc$.

Let us mention that, although we require that $\cc$ is real, the
results in this paper are also valid for non-real plane algebraic
curves. In addition, the condition  $(1:0:0),(0:1:0)\not\in \cc^h$
can be avoided by performing a suitable affine orthogonal linear
change of coordinates. The requirement on the smoothness of
$\cc^h$ at infinity, might be avoided by performing a suitable
projective linear change of coordinates. However, differently to
affine orthogonal linear changes,  in general,  projective changes
of coordinates do not preserve properly the closeness between the
input and output curves.

\section{Preliminaries on $\epsilon$-points}\label{section-e-points}

Our fundamental technique to deal with the approximate
parametrization problem is the use of $\epsilon$-points. The
notion of $\epsilon$--point of an algebraic variety  was
introduced by the authors (see \cite{PSS}, \cite{PSS2},
\cite{PSS3}) as a generalization of the notion of approximate root
of a univariate polynomial.  In this section, we briefly summarize
some previous notions introduced in \cite{PSS} and \cite{PSS2},
and geometric properties obtained in \cite{PSS3}. We start with
the notion of proper degree.

\parasmall

\begin{definition}\label{def-proper-degree}  We say that a polynomial
$g\in {{\Bbb C}}[x,y]$ has {\sf proper degree}  $\ell$ if the total
degree of $g$ is $\ell$,  and $\exists$ $\overrightarrow{v} \in {{\Bbb
N}}^2$, with $|\overrightarrow{v}|=\ell$, such that $
|g^{\overrightarrow{v}}|
>\epsilon \|g\|.$

 We say that an algebraic plane curve has {\sf proper degree} $\ell$ if its
defining polynomial has proper degree $\ell$. \qed
\end{definition}

\parasmall

The notion of $\epsilon$--point is as follows.

\parasmall

\begin{definition}\label{def-e-point}
${P}  \in {{\Bbb C}}^2$ is an {\sf $\epsilon$--(affine)
point} of  $\cal C$ if $\left|f(P )\right|<\epsilon \|f\|$. \qed
\end{definition}

\parasmall

 In this situation, we   introduce the notion
 of  $\epsilon$-singularity, pure $\epsilon$-singularity, and $\epsilon$-ramification point.

\parasmall

\begin{definition} \label{def--mult-e-point}  Let  ${P}  \in {{\Bbb
C}}^2$ be an $\epsilon$--point of $\cal C$.
\begin{itemize}
\item[(i)] The
{\sf $\epsilon$-multiplicity} of $P $ on $\cc$ (we denote it by $\emul(P ,\cc)$) is the smallest natural number  $r
\in {{\Bbb N}}$ satisfying that
\begin{itemize}
\item[(1)]  $\forall\, \overrightarrow{v} \in {{\Bbb
N}}^2$, such that $0\leq |\overrightarrow{v}| \leq r-1$, it holds
that $ |f^{\overrightarrow{v}}(P )|<\epsilon \|f\|$,
\item[(2)]  $\exists\, \overrightarrow{v} \in {{\Bbb
N}}^2$, with $|\overrightarrow{v}|=r$, such that $
|f^{\overrightarrow{v}}(P )|\geq \epsilon \|f\|. $
\end{itemize}
\item[(ii)]    $P $ is an {\sf $\epsilon$--(affine) simple point}
of $\cal C$ if $\emul(P ,\cc)=1$; otherwise,  $P $ is an  {\sf
$\epsilon$--(affine) singularity} of $\cal C$.
\item[(iii)] ${P} $ is a {\sf $k$-pure $\epsilon$--singularity of
 $\cal C$}, with $k\in\{1,2\}$,  if $\emul(P ,\cc)>1$ and $|f^{\emul(P ,\cc)\cdot
\overrightarrow{e_k}}(P )| \geq \epsilon \|f\|.$
\item[(iv)] ${P} $ is an {\sf $\epsilon$--(affine)
ramification point} of  $\cal C$  if $\emul(P ,\cc)=1$, and
either $|f^{\overrightarrow{e_1}}(P )| < \epsilon \|f\|$ or
$|f^{\overrightarrow{e_2}}(P)| < \epsilon \|f\|.$
\qed
\end{itemize}
\end{definition}

\parasmall

Note that, since $\cc$ has proper degree,  $0\leq \mul(P ,\cc)\leq \emul(P ,\cc) \leq \deg(\cc)$, where
 $\mul(P ,\cc)$ denotes the ``exact" multiplicity of $P $ on $\cc$.
 For instance,  the origin has exact multiplicity 1, and $\epsilon$-multiplicity 2,
 on the curve defined by $\frac{\epsilon}{2} x+x^3+y^2$.
 In the exact case, if $\cc$ is irreducible, $\mul(P ,\cc)<\deg(\cc)$.
 Thus one may expect that in the approximate case, if $\cc$ is $\epsilon$-irreducible,
 then $\emul(P ,\cc)<\deg(\cc)$. Although this is the case in all the examples we have tried,
 we have not been able to prove it. So in this paper, when computing $\epsilon$-multiplicities,
 we also consider the possibility $\emul(P ,\cc)=\deg(\cc)$.

\parasmall

The following lemma is a direct generalization of Lemma 3 in \cite{PSS}.

\parasmall

\begin{lemma}\label{lema-discos-multiplicidades}
Let  ${P}  \in {{\Bbb
C}}^2$ be an $\epsilon$--point of $\cal C$. There exists $\delta>0$ such that every ${Q} \in {{\Bbb
C}}^2$, satisfying that $\|P -Q\|_2<\delta$, is an $\epsilon$-point of $\cc$ with $\emul(Q,\cc)\geq\emul(P ,\cc)$.
\end{lemma}

\noindent {\bf Proof.} Simply observe that the reasoning of Lemma
3 in \cite{PSS} is also valid over $\Bbb C$. \qed

 \parasmall

The following example shows that, in Lemma \ref{lema-discos-multiplicidades},  the $\epsilon$-multiplicity of $Q$ can be strictly bigger than $\emul(P ,\cc)$.

\parasmall

\begin{example}\label{Ejemplo-multi-degenerado}
Let $\cc$ be defined by $f(x,y)=x^3y+y^3x+x^3+\frac{\epsilon}{2}x^2+\epsilon y+\frac{\epsilon}{2}$; note that $\|f\|=1$. For $P =(0,0)$, one has $$f(P )=\frac{\epsilon}{2},\quad f^{\overrightarrow{e_1}}(P )=0,\quad f^{\overrightarrow{e_2}}(P )=\epsilon.$$ So, $\emul(P ,\cc)=1$. Now, we consider the sequence of points $\{P_n=(-\frac{1}{n},0)\}_{n\geq 1}$. Then,
  $$f(P_n)=\frac{\epsilon}{2}+\frac{\epsilon}{2n^2}-\frac{1}{n^3},\quad f^{\overrightarrow{e_1}}(P_n)=\frac{3}{n^2}-\frac{\epsilon}{n}, \quad f^{\overrightarrow{e_2}}(P_n)=\epsilon-\frac{1}{n^3},$$  $$ f^{\overrightarrow{(2,0)}}(P_n)=\epsilon-\frac{6}{n}, \quad f^{\overrightarrow{(1,1)}}(P_n)=\frac{3}{n^2},\quad f^{\overrightarrow{(0,2)}}(P_n)=0,\quad f^{\overrightarrow{(3,0)}}(P_n)=6$$
So, for $n$ sufficiently large, $\emul(P_n,\cc)=3$. \qed
\end{example}

\parasmall

\begin{definition}\label{def-depth}
  Let $P $ be an
$\epsilon$-point of $\cal C$ and $r=\emul(P ,\cc)$. If  $P $ is
$k$--pure, with $k\in\{1,2\}$, we  define the {\sf $k$-weight of
$P $} as
\[\we_k(P)=\displaystyle{\max_{i=0,\ldots,r-1}\left \{
\left|\frac{r!\cdot f^{i\cdot \overrightarrow{e_k}}(P )}{i!\cdot
f^{r\cdot \overrightarrow{e_k}}(P )}\right|^{\frac{1}{r-i}}\right
\}}. \] We define the {\sf weight of $P $}, denoted by $\we(P )$,
as $\max\{\we_1(P ), \we_2((P )\}$, if $P$ is pure in both
directions,  and as the corresponding $k$-weight otherwise. \qed
\end{definition}

\parasmall

The following two rational functions were introduced in \cite{sa:te}, and will
 play an important role in this development:
\[{\cal R}_{\inn}(x)=2x
\left(\frac{1}{1+3x}+\frac{16x}{(1+3x)^3}\right),\hspace*{7mm} {\cal
R}_{\out}(x)=
\frac{1}{2}-\frac{x(1-9x)}{2(1+3x)}-\frac{32x^2}{(1+3x)^3}.\]
Furthermore, these two rational functions give information on how close an $\epsilon$--point is
 to an exact point of the curve $\cal C$ (see Theorem 2 and Corollary 3 in \cite{PSS3}).

\section{Preliminaries on Symbolic Parametrization}\label{section-symbolic-parametrization}

In this section, we briefly recall the symbolic parametrization algorithm for rational plane algebraic curves of degree $d>2$ (note that lines and conics can be trivially parametrized by lines) based on $(d-2)$ adjoint curves; for further details see \cite{SWP}. For this purpose, throughout this section we assume that $\cc$ is rational (i.e. its genus is zero).
In addition,  taking into account our requirements in Section \ref{seccion-approximate-parametrization} and for simplicity sake, we   assume in this section that all singularities are affine and ordinary. Again, for a complete description see \cite{SWP}.

\parasmall

The idea is to use a linear system of curves  such that for almost
every curve in this  system, all its intersections with $\cc^h$,
except one, are predetermined; recall that $\cc^h$ is the
projective closure of $\cc$. Moreover, the set of all these
intersection points is the same one for every curve in the system,
and the points in this set are called the ``base points". Thus, if
one computes the intersection points of $\cc^h$ with a generic
representative of the system, the expression of the unknown
intersection point gives the parametrization of the curve in terms
of the parameter defining the linear system.

\parasmall

More precisely, let ${\cal H}_{d-2}$ be the linear
system of adjoint curves to $\cc^h$ of degree $d-2$. That is,
${\cal H}_{d-2}$ is the linear system of  curves of degree $d-2$
having each $r$--fold of $\cc^h$ as a base point of multiplicity  $r-1$; i.e. as a point of multiplicity at least $r-1$. In particular it implies that the
multiplicity of intersection of a curve in ${\cal H}_{d-2}$ and
$\cc^h$ at a base point of multiplicity $r-1$ is at least $r(r-1)$.
Thus, using that the genus of $\cal C$ is zero, and  taking into account B\'ezout's
Theorem, one deduces that $d-2$ intersections of $\cc^h$ and a generic element in ${\cal H}_{d-2}$ are not predetermined.
 In this situation, one may take $(d-3)$ simple points on $\cc^h$, and determine the 1-dimensional linear
subsystem  ${\cal H}_{d-2}^*$ of ${\cal H}_{d-2}$ obtained when these simple points are  required to be base points of multiplicity 1.
In this way,  the number of predetermined intersections (counted with multiplicity) is $(d-1)(d-2)+(d-3)$, i.e. only one intersection point is missing. Thus, computing this
free intersection one finds a rational parametrization of $\cc^h$. Summarizing these ideas one has the following
process:


\begin{itemize}
\item[(1)] Compute the  singularities of $\cc^h$ as well as  their multiplicities
 (recall that we have assumed that all singularities are affine and ordinary).
\item[(2)] Determine the linear system ${\cal H}_{d-2}$ of adjoint curves of degree $(d-2)$ to
$\cc^h$.
\item[(3)] Compute  $d-3$  different simple  points on ${\cal C}^h$.
\item[(4)] Determine the linear subsystem   ${\cal H}^{*}_{d-2}$ of ${\cal H}_{d-2}$ by requiring that every simple points in Step 3 is a base
point of multiplicity one.
\item[(5)] Compute the free intersection point of ${\cal H}^{*}_{d-2}$ and
$\cc^h$.
\end{itemize}


Let us make a comment on how to computationally perform the steps
in the above process. Step (1) can be performed, for instance,
using resultants. In Step (2), one considers a homogeneous
polynomial $H(x,y,z)$ of degree $(d-2)$ with undetermined
coefficients. Now, for each singular point $P$ of multiplicity $r$
one requires that $H$ and all its partial derivatives till order
$(r-1)$ vanish at $P$. This generates a linear system of equations
in the undetermined coefficients of $H$. Solving it, and
substituting in $H$, we get the defining polynomial of  ${\cal
H}_{d-2}$; let us call it again $H$. Step (3) may be performed by
intersecting $\cc^h$ with lines (see \cite{SWP} for advanced
approaches); although it is not necessary, looking for the
parallelism with the reasoning in Section
\ref{seccion-approximate-parametrization}, we take  affine simple
points. Step (4) can be approached as Step (2), i.e. requiring
that $H$ vanishes at each simple point, solving the provided
linear system and substituting the solution in $H$; let
$H^*(t,x,y,z)$ be the defining polynomial of ${\cal H}_{d-2}^*$
(note that $\dim({\cal H}_{d-2}^*)=1$) . Finally, let us deal with
Step (5). For this purpose, let
$\{Q_i:=(q_{i,1}:q_{i,2}:1)\}_{i=1,\ldots,s}$ be the singularities
and $r_i$ the multiplicity of $Q_i$. Also, let
$\{P_i:=(p_{i,1}:p_{i,2}:1)\}_{1,\ldots,d-3}$ be the simple points
determined in Step (3). Then,  the free intersection point is
obtained by computing the primitive part, w.r.t. $t$, of the
resultants of $H^*(t,x,y,1)$ and ${f}(x,y)$ with respect to $x$
and $y$, respectively. Indeed, it holds that  (see  \cite{SSS05})
\[\begin{array}{l}
S_1(x,t)=\resultant_y(H^*(t,x,y,1), {f}(x,y))=\prod_{i=1}^{s}(x-{q}_{i,1})^{r_i(r_i-1)}
\prod_{i=1}^{d-3}(x-p_{i,1})M_1(x,t), \\ S_2(y,t)=\resultant_x(H^*(t,x,y,1), {f}(x,y))=\prod_{i=1}^{s}(y-{q}_{i,2})^{r_i(r_i-1)}
\prod_{i=1}^{d-3}(y-{p}_{i,2})M_2(y,t),
\end{array} \] where
$\degree_x(M_1)=\degree_y(M_2)=1$. Therefore, the parametrization is the solution in $\{x,y\}$ of   $\{M_1(x,t)=0,  M_2(y,t)=0\}$.

\section{$\epsilon$-Rational Curves}\label{seccion-e-racional}

In this section we introduce the notion of $\epsilon$-rationality
of a plane algebraic curve. This notion plays the corresponding
role in the approximate frame that the rationality does for exact
algebraic curves. We will have two main difficulties.  On one
hand, computing the $\epsilon$-multiplicity and on the other,
differently to the exact case, we will have in general more
$\epsilon$-singularities than expected, and we will need to
associate them; we will solve this last difficulty introducing a
suitable concept of cluster.

\parasmall

We first need to determine the $\epsilon$-singularities. To check
the existence and perform the actual computation of the
$\epsilon$-singularities, one has to solve the system of algebraic
equations
$$\{f^{\overrightarrow{0}}(x,y)=0,\quad f^{\overrightarrow{e_1}}(x,y)=0,\quad f^{\overrightarrow{e_2}}(x,y)=0\},$$
 under
fixed precision $\epsilon \|f\|$. This can be done, for instance, by applying root
finding techniques (see \cite{cor:gi:tra:wat}, \cite{go:van}, \cite{Ho}).
Note that since $f$ is $\epsilon$-irreducible then it is irreducible, and hence
the above system has finitely many solutions. Let ${\cal S}_1$ be
the set of solutions. One may accelerate the computation by
working (if possible) with two co-prime polynomials, instead of
three, to get a finite super-set of the set of solutions, from
where the $\epsilon$-singularities  are detected  afterwards.

\parasmall

Now, for $P\in {\cal S}_1$, we want to compute $\emul(P,\cc)$.
This can be obviously done by substituting $P$ at the
corresponding partial derivatives and checking the conditions in
Def. \ref{def--mult-e-point} (1). Seemingly, there is no
difficulty on that. However, in Lemma
\ref{lema-discos-multiplicidades}, we have seen that for each
$\epsilon$-point $P$ of $\epsilon$-multiplicity $r$ there exists
an open disk $U$ centered at $P$  such that if $Q\in U$, then $Q$
is an $\epsilon$-point of $\epsilon$--multiplicity at least $r$.
So, an small perturbation of $P$ may produce an incorrect answer
for the $\epsilon$-multiplicity; see, for instance, Example
\ref{Ejemplo-multi-degenerado}. We are indeed interested in
assigning the maximum possible $\epsilon$-multiplicity to the
$\epsilon$-point. The proof of Lemma 3 in \cite{PSS}, and hence of
Lemma \ref{lema-discos-multiplicidades}, shows how to detect the
radius of one of these open disks, so one may try to estimate the
maximum $\epsilon$-multiplicity at the disk. Nevertheless, in
practice, this is unfeasible. Instead, we propose a different
strategy that, although it does not ensure the achievement of the
maximum,  in practical examples turns to work efficiently.

\parasmall

More precisely, for each $k\in \{2,\ldots,d-1\}$, we take
$\overrightarrow{u_1},\ldots,\overrightarrow{u_s}\in {\Bbb N}^2,$ with $2\leq s \leq k+1$ (in practice $s=2$) such that for all $i$, $|\overrightarrow{u_i}|=k$ and $\gcd(f^{\overrightarrow{u_1}},\ldots,f^{\overrightarrow{u_s}})=1$, and we solve
$\{f^{\overrightarrow{u_1}}=0,\ldots,f^{\overrightarrow{u_s}}=0\},$
 under fixed precision $\epsilon \|f\|$. Let ${\cal A}_k$ be the set of  solutions. Then,  for $k\in \{2,\ldots,d-1\}$ we consider the set (note that ${\cal S}_1$ is defined above)
  \[ {\cal S}_k =\{  P\in {\cal A}_k\,/\, |f^{\overrightarrow{w}}(P)|< \epsilon \|f\|\,\,\,\forall \overrightarrow{w} \in {\Bbb N}^2\,\,\mbox{with}\,\,  |\overrightarrow{w}|\leq k\}. \]
If for a given $k$ and for all $s$ it holds that $\gcd(f^{\overrightarrow{u_1}},\ldots,f^{\overrightarrow{u_s}})\neq 1$,   we take ${\cal S}_k=\emptyset$. Finally we consider the set
\[ {\cal S}=\bigcup_{k=1}^{d-1} {\cal S}_k. \]
It is clear that in general we introduce additional $\epsilon$-singularities, and we will have to generate a process (the cluster construction) to identify
them. Nevertheless, each new $\epsilon$-singularity, after identification, will increase the $\epsilon$-multiplicity of the original one.

\parasmall

\begin{definition}\label{Def-e-singular-locus}
The  set  ${\cal S}$, introduced above, is called the {\sf
$\epsilon$-(affine)-singular locus of $\cal C$}. We denote it by
$\sing({\cal C})$. \qed
\end{definition}

\parasmall

\begin{example}\label{Ejemplo-multi-degenerado-segunda-parte}
Let us take $\epsilon=0.001$ in Example \ref{Ejemplo-multi-degenerado}.
The $\epsilon$-singular locus of this curve is $\sing(\cc)={\cal S}_1\cup{\cal S}_2\cup {\cal S}_3$ where
\[
\begin{array}{l}
\begin{array}{ll}
{\cal S}_1=\{\hspace*{-0.3cm}&P_1=(0.02131893405+0.009609927603i, 0.02442855631+0.1171004584i),\\
          \,&P_2=(0.004713033954+0.02355323617i, -0.07491796596-0.09032199938i),\\
           \,&P_3=(-0.01424770212+0.01818884517i, 0.1084633939+0.05315246871i),\\
           \,& P_4=(-0.02443272919, -0.1159479025),\\
\,&P_5=(-0.01424770212-0.01818884517i, 0.1084633939-0.05315246871 i),\\
 \,&P_6=(0.004713033954-0.02355323617 i, -0.07491796596+0.09032199938 i),\\
 \, &P_7=(0.02131893405-0.009609927603 i, 0.02442855631-0.1171004584 i)\},
\end{array}
\\
\begin{array}{ll}
{\cal S}_2=\{\hspace*{-0.3cm}&P_8=(-0.0001666666667, 0)\}.
\end{array}
\\
\begin{array}{ll}
{\cal S}_3=\emptyset.
\end{array}
\end{array}
\]
Moreover, $\emul(P_1)=\cdots=\emul(P_7)=2$ but $\emul(P_8)=3$. Note that considering only ${\cal S}_1$ we would have not found a point with $\epsilon$-multiplicity $3$.
\qed
\end{example}

\parasmall

As we could check in the previous example,  the difficulty appears
when observing that we may have two (in general more than two)
$\epsilon$-singularities $P$ and $Q$ that are very ``{\it close}",
and somehow we need  to identify them. To approach this, we
introduce the notion of cluster of $\epsilon$-singularities.
Intuitively,  two $\epsilon$-singularities $P$ and $Q$ of $\cal C$
are in the same cluster, if the disks centered at $P$ and $Q$ (of
certain radius)  are a small vibration of each other. The radius
and the vibration are measured by means of the value of the
function ${\cal R}_{\out}$  at the weight and the tolerance,
respectively (see Section \ref{section-e-points}). Since the
notion of weight requires that the $\epsilon$-singularities are
pure, for non-pure $\epsilon$-singularities we will take radius
zero. More precisely, we introduce the following definition.

\parasmall

\begin{definition}\label{def--radius}
Let $P$ be an $\epsilon$-point of $\cal C$. We define its {\sf radius},
and we denote it by $\rad(P)$, as ${\cal R}_\out(\we(P))$ if $P$ is pure and zero otherwise. \qed
\end{definition}

\parasmall

\begin{definition}\label{def--cluster} Let ${\cal A}$ be a finite set of ${\epsilon}$--points  of ${\cal
C}$. For  $P\in {\cal A}$  we define the {\sf cluster of $P$ w.r.t. $\cal A$} as the set of all points $Q\in \cal A$
such that  at least one
of the following conditions is verified:
\begin{itemize}
\item[(1)] $\|P-Q\|_2 {+}\mid \rad(P)-\rad(Q) \mid< {\cal R}_{\out} (\epsilon ),$
\item[(2)]  there exists  $P'\in \cal A$ such that
$\|P'-P\|_2 {+}|\rad(P')-\rad(P)|< {\cal R}_{\out} (\epsilon )$
 and
 $\|P'-Q\|_2{+}|\rad(P')-\rad(Q)|< {\cal R}_{\out} (\epsilon ).$
\end{itemize}

We say that $R$ is  {\sf a candidate to be the representative of a
cluster}, if $R$ is a point of the cluster of maximum
$\epsilon$--multiplicity. We say that $R$ is  {\sf a
representative of a cluster} if it is a candidate and  $|f(R)|
\leq | f(Q) |$ for all the other candidates $Q$. We define the
{\sf $\epsilon$--multiplicity of the cluster} as the
$\epsilon$--multiplicity of any of its representatives.

We denote a
cluster by $\cluster_r(R,{\cal A})$, where $r$  is the $\epsilon$--multiplicity and $R$ a
representative, and by $\cluster_r(R)$ when ${\cal A}=\sing(\cc).$ \qed
\end{definition}

\parasmall

Now, we are ready to introduce the notion of $\epsilon$-rationality.

\parasmall

\begin{definition}\label{D-erational}
If $\{\cluster_{r_i}(P_i)\}_{i=1,\ldots,s}$ is the cluster
decomposition of $\sing(\cc)$,  we say
that $\cal C$ is  {\sf $\epsilon$-(affine) rational} if $
(d-1)(d-2)-\sum_{i=1}^{s} r_i(r_i-1) =0. $\qed
\end{definition}

\parasmall

\begin{remark}
Note that in the previous theoretical development we have not considered singularities (neither $\epsilon$-singularities) at infinity. We leave this extension of the concept of $\epsilon$-rationality for further research.\hfill \qed
\end{remark}

\parasmall

If we apply the previous ideas to Example
\ref{Ejemplo-multi-degenerado-segunda-parte} (see also Example
\ref{Ejemplo-multi-degenerado}), with $\epsilon=0.001$ we get that
the 8 points of $\sing(\cc)$ belong to the same cluster. So, the
cluster decomposition is
$\{\cluster_{3}(P_8)=\{P_1,\ldots,P_8\}\}$. Therefore, $\cc$ is
$\epsilon$-rational; indeed, it is $\epsilon$-monomial, and thus
parametrizable with the techniques in \cite{PSS}. We finish the
section with a more general example.

\parasmall

\begin{example}
Let us consider $\epsilon= 0.005$ and the curve $\cc$ of proper
degree $5$ defined by the polynomial (see Fig.\ref{clusters-primera-parte}):
 \vspace{1 mm}

\parbox{15cm}{\small
$f(x,y)=-2.199771784x^2-0.2197717843x^4y-0.9016804979x^3y^2+1.858817427x^3-1.891680498y^4
+0.9899999999xy^3+0.9899999999x^2y-1.055726141x^2y^3+0.3409543568y^2+0.9899999999x^4
+0.9899999999xy^4+0.9899999999y^3-0.1869087137x^5+5.235497925xy^2-1.770497925x^2y^2
+1.45213693x^3y-0.1440456432xy-0.52786307y^5+0.01
 $.}

\vspace{1 mm}

\noindent
The $\epsilon$-singular locus is $\sing(\cc)={\cal S}_1\cup {\cal S}_2\cup {\cal S}_3$, where
\[
\begin{array}{l}
\begin{array}{ll}
{\cal S}_1=\{\hspace*{-0.3cm}&P_1=(-0.9956027274+0.0004067223817 i,
0.001447687187+0.9982777543 i),\\
          \,&P_{2}=(1.011706789-0.1320874194 i, -1.008532436+0.06832949372 i),\\
           \,& P_3=(1.007458642, -1.044045331),\, P_4=(0.9909273695, -0.9540334161) ,\\
\,&P_5=(1.011706789+0.1320874194 i, -1.008532436-0.06832949372 i),\\
 \,&P_6=(-0.9956027274-0.0004067223817 i, 0.001447687187-0.9982777543 i),\\
 \, &P_7=(0, 0) \},
\end{array}
\\
\begin{array}{ll}
{\cal S}_2=\{\hspace*{-0.3cm}&P_{8}=(1.000000001, \,-1.) \},
\end{array}
\\
\begin{array}{ll}
{\cal S}_3=\emptyset.
\end{array}
\end{array}
\]
Moreover, $\emul(P_1)=\emul(P_2)=\emul(P_7)=2$, and $\emul(P_3)=\emul(P_4)=\emul(P_5)=\emul(P_6)=\emul(P_8)=3$. Furthermore, the
cluster decomposition is  (see Fig. \ref{clusters-primera-parte}):
\[\begin{array}{l}
\cluster_2(P_{1}) = \{P_1\}, \\ \cluster_2(P_{2})=\{P_2\}, \\
\cluster_2(P_{7})= \{P_7\}, \\
\cluster_3(P_8) = \{P_3, P_4, P_5, P_6, P_8\}.
\end{array} \]
 Thus, $\cal C$ is
$\epsilon$-rational. 
\begin{center}
\begin{figure}[ht]
\centerline{
\psfig{figure=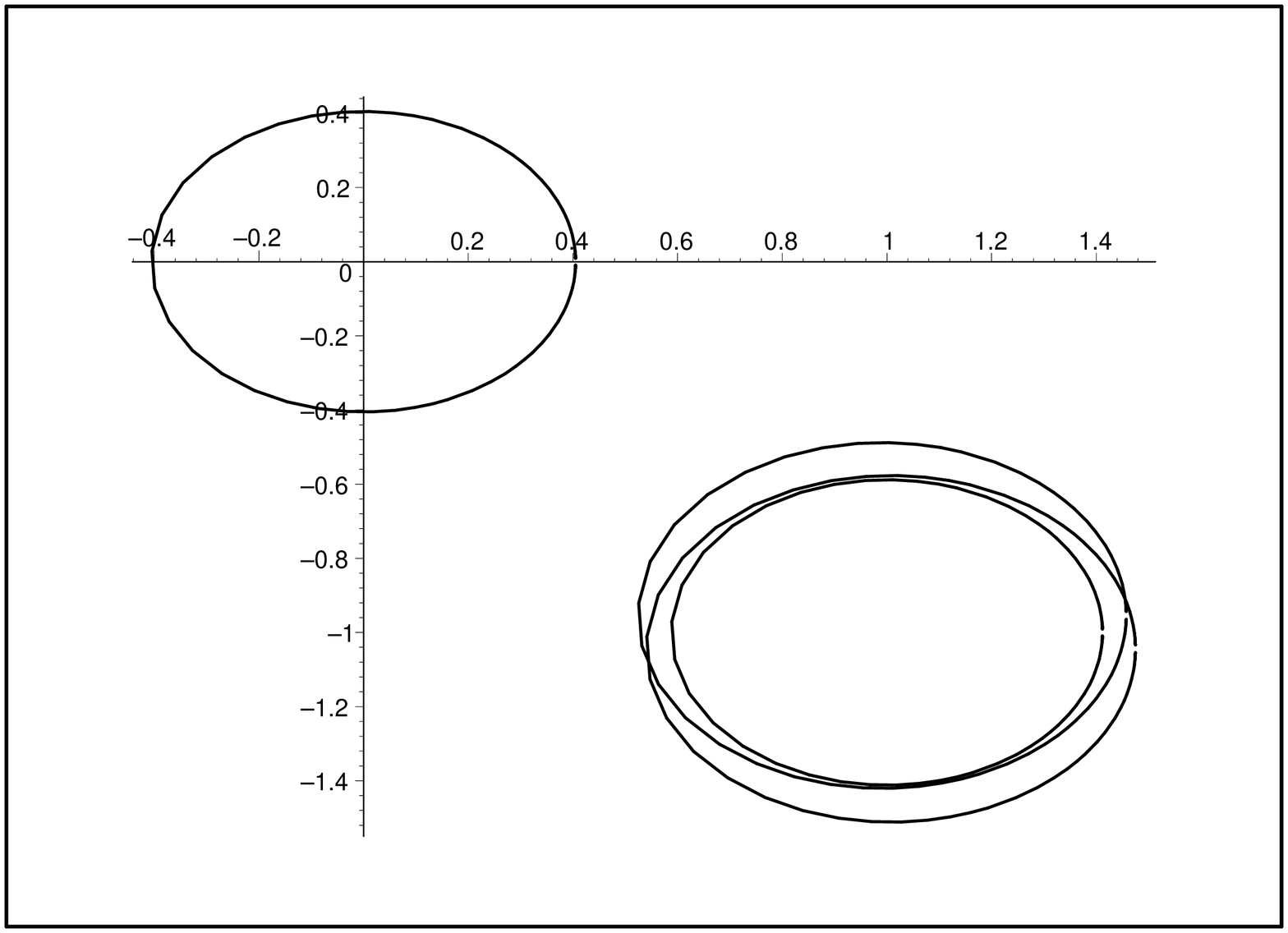,width=4.4cm,height=4.4cm,angle=270}
\psfig{figure=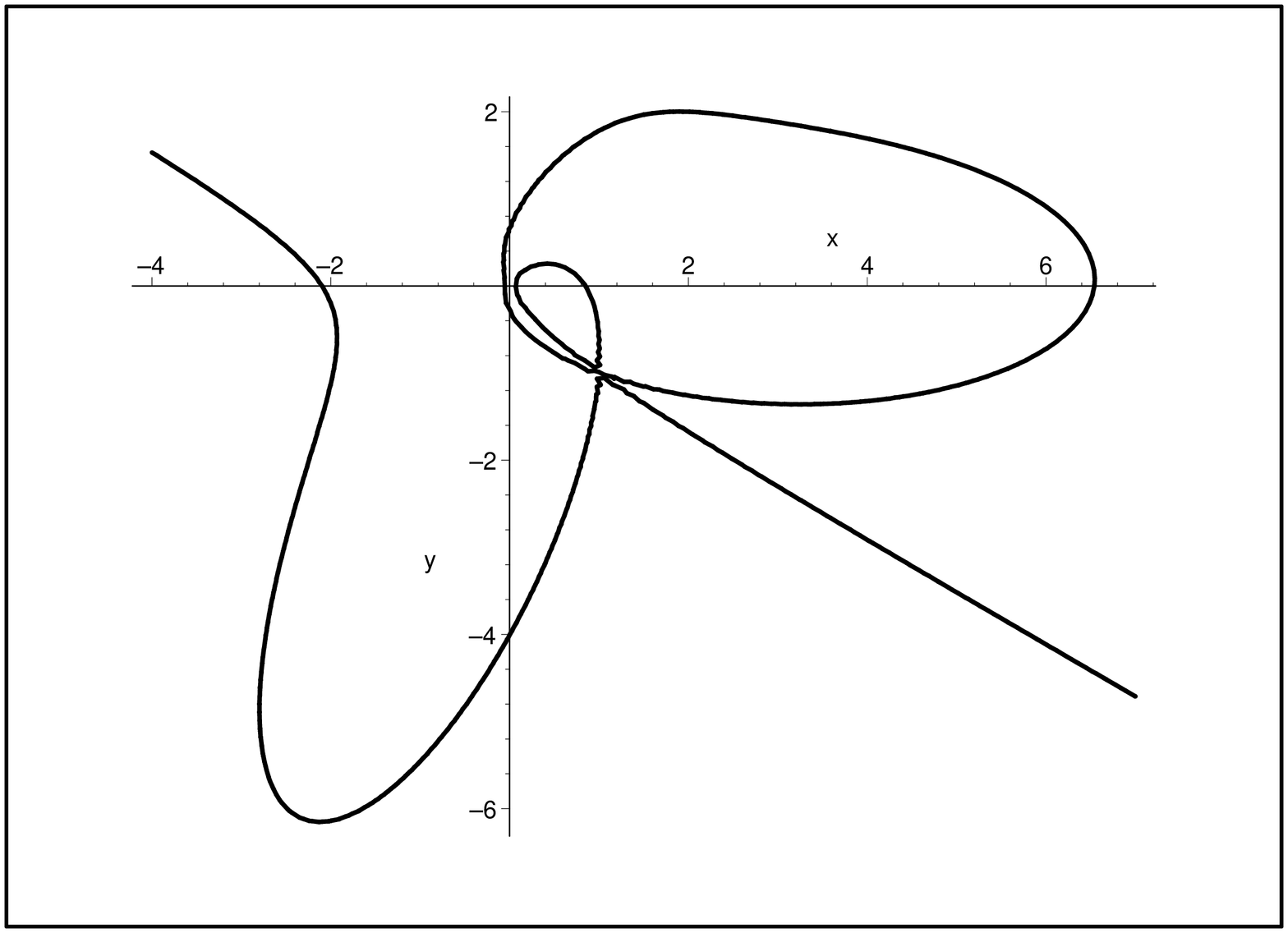,width=4.4cm,height=4.4cm,angle=270} }
\caption{{\sf Left:} Clusters. {\sf Right:} Curve
$\cc$}\label{clusters-primera-parte}
\end{figure}
\end{center}\vspace*{-1cm}\hfill \qed
\end{example}

\section{Approximate Parametrization Algorithm}\label{seccion-approximate-parametrization}

In this section, we present our approximate parametrization algorithm. For this purpose, we assume that $\cc$ is $\epsilon$-rational of proper degree $d>2$ (note that for $d=1$ the problem is trivial, and for $d=2$ one can apply the algorithm in \cite{PSS}), and  that
$$\{\cluster_{r_i}(Q_i)\}_{i=1,\ldots,s},\,\,\mbox{where}\,\,Q_i:=(q_{i,1}:q_{i,2}:1),$$
is the cluster decomposition of $\sing(\cc)$. Furthermore, if
possible, i.e. when there exists a real representative of the
cluster, we take $Q_i$ real.

\parasmall

In this situation, we adapt the algorithm in Section \ref{section-symbolic-parametrization} as follows.
 Let $\cc^h$ be the projective closure of $\cc$. We consider the linear
system of   curves $\overline{\cal H}_{d-2}$ of degree
$(d-2)$ given by the divisor $\sum_{i=1}^{s} r_i Q_i$. That is,  ${Q}_i$ is a base point of (exact) multiplicity
$r_i-1$ of the linear system. Afterwards, one computes  $(d-3)$ $\epsilon$--simple affine points on $\cc^h$ (see below for details), and  determines the linear subsystem $\overline{\cal
H}^{\,*}_{d-2}$ of $\overline{\cal H}_{d-2}$ obtained by intersecting $\overline{\cal H}_{d-2}$ with the linear system of $(d-2)$-degree curves generated by   the divisor $\sum_{i=1}^{d-3}P_i$; say that $P_i:=(p_{i,1}:p_{i,2}:1)$. If ${P}_i,\,{Q}_j$ would be exact points and singularities,
respectively, of  $\cc^h$, then $\dim(\overline{\cal
H}^{\,*}_{d-2})=1$ (see Chap. 4 in \cite{SWP}). However, in our case, since we are working  with $\epsilon$-points we can only ensure that $\dim(\overline{\cal
H}^{\,*}_{d-2})\geq 1$ (see Theorem 2.56 in \cite{SWP}). If this dimension is strictly bigger than 1, we can either take more $\epsilon$-simple points till dimension 1 is reached, or we can take an small perturbation of the $\epsilon$-points such that the effective divisor $\sum_{i=1}^{s} r_i Q_i  + \sum_{i=1}^{d-3}P_i$ is in general position (see page 49 in \cite{SWP}), and hence the dimension is 1. So, we can assume w.l.o.g. that $\dim(\overline{\cal
H}^{\,*}_{d-2})=1$.  Let, then,  $\overline{H}_{d-2}(t,x,y,z)$ be the defining homogeneous polynomial of  $\overline{\cal H}^{*}_{d-2}$.

\parasmall

 At this point,
if ${P}_i,\,{Q}_j$ would be exact points and singularities,
respectively, of  $\cal C$, the symbolic algorithm
presented in Section \ref{section-symbolic-parametrization} would output the parametrization
${\cal P}(t)=(\frac{p_1(t)}{q_1(t)},\frac{p_2(t)}{q_2(t)})$, where
\[ q_1(t)x-p_1(t)= \frac{\resultant_y(\overline{H}^{*}_{d-2}(t,x,y,1), {f}(x,y)) }{\prod_{i=1}^{s}(x-{q}_{i,1})^{r_i(r_i-1)}\prod_{i=1}^{d-3}(x-{p}_{i,1})}, \]
\[ q_2(t)y-p_2(t)= \frac{\resultant_x(\overline{H}^{*}_{d-2}(t,x,y,1), {f}(x,y)) }{\prod_{i=1}^{s}(y-{q}_{i,2})^{r_i(r_i-1)}\prod_{i=1}^{d-3}(y-{p}_{i,2})}.  \]
However, in our case, ${P}_i,\,{Q}_j$ are not exact points, but
$\epsilon$--points. So these rational functions are not, in
general, polynomials. Nevertheless, considering if necessary a
small perturbation of $\overline{H}_{d-2}$, the quotient of the
division of each numerator by its denominator is linear as
polynomial in either $x$ or $y$. Then, the idea is to determine
the parametrization from these linear quotients. For this purpose,
we will consider (if necessary) two  perturbations, both affecting
$\overline{H}_{d-2}$. The first one will ensure that the degree in
the resultants is the expected one, namely $d(d-2)$. The second
will guarantee that the output is indeed a parametrization; i.e.
that not both components are constants. Note that, in the exact
case, these two facts are provided  by the theory.

\parasmall

More precisely, let $\overline{H}^{*}_{d-2}(t,x,y,z)=H_1(x,y,z)+t
\, H_2(x,y,z)$, and let ${\cal D}_i$ be the projective curve
defined by $H_i$, $i=1,2$. We recall that $(1:0:0), (0:1:0)\not\in
\cc^h$. Now, we need to ensure that either $\cc^h, {\cal D}_1$  or
$\cc^h, {\cal D}_2$ do not have common points at infinity. If this
is not the case, let $\{R_1,\dots,R_m\}$ be the points of $\cc$ at
infinity and $K(\rho_1,\rho_2,x,y,z)=\rho_1 x^{d-2}+\rho_2
y^{d-2}$, where $\rho_i$ are parameters. Then, we consider in
${\Bbb C}^2$ the union $\cal L$ of the affine lines defined by
$H_2(R_i)+K(\rho_1,\rho_2,R_i)=0$, for $i=1,\ldots,m$. Note that,
since $R_i$ are points at infinity, the polynomials
$H_2(R_i)+K(\rho_1,\rho_2,R_i)\in {\Bbb C}[\rho_1,\rho_2]$ are not
constant, and hence define lines. So, taking values for
$\rho_1,\rho_2$ (say, small real numbers) we consider an small
perturbation that ensures that the above requirement is satisfied.

\parasmall

Thus, in what follows we assume that ${\cal D}_2$ and $\cc^h$ do
not have common points at infinity. Therefore, if $F$ is the
homogenization of $f$,    by Lemma 3.1 in \cite{AS}, one has that
\[ \deg_x(\resultant_y(\overline{H}^{*}_{d-2}, {F}))=\deg_y(\resultant_x(\overline{H}^{*}_{d-2}, {F}))=d(d-2). \]
Moreover, since $\overline{\cal
H}^{\,*}_{d-2}$ and ${\cal C}^h$ do not have common points at infinity, it holds that
\[ \deg_x(\resultant_y(\overline{H}^{*}_{d-2}(t,x,y,1), {f}))=\deg_y(\resultant_x(\overline{H}^{*}_{d-2}(t,x,y,1), {f}))=d(d-2). \]
Now, we consider the polynomials
$$A_1(x)=\prod_{i=1}^{s}(x-{q}_{i,1})^{r_i(r_i-1)} \prod_{i=1}^{d-3}(x-{p}_{i,1}),\, A_2(y)=\prod_{i=1}^{s}(y-{q}_{i,2})^{r_i(r_i-1)}\prod_{i=1}^{d-3}(y-{p}_{i,2}). $$
Since $\cc$ is $\epsilon$-rational, it holds that
\[ \deg_x(A_1(x))=
\deg_y(A_2(y))=d(d-2)-1. \]
Let   $B_1(x,t):=\overline{q}_1(t)x-\overline{p}_1(t)$  be the quotient of $S_1(x,t):=\resultant_y(\overline{H}^{*}_{d-2}(t,x,y,1), {f}(x,y))$ and $A_1(x)$. Similarly let  $B_2(y,t):=\overline{q}_2(t)x-\overline{p}_2(t)$ be the quotient of
$S_2(y,t):=\resultant_x(\overline{H}^{*}_{d-2}(t,x,y,1), {f}(x,y))$  and $A_2(y)$. Then, we output
$$\overline{{\cal P}}(t)=\left(\frac{\overline{p}_1(t)}{\overline{q}_1(t)},\frac{\overline{p}_2(t)}{\overline{q}_2(t)}\right)$$
as approximate parametrization of $\cc$.

\parasmall

Intuitively one sees that, in practice,  $\overline{\cP}(t)$ will
be always a parametrization. In order to prove this claim, we
repeat the reasoning but introducing a new perturbation of
$\overline{\cal H}^{*}_{d-2}$. More precisely, let
$\Delta=(\delta_1,\ldots,\delta_6)$ be a family of perturbing
parameters and  let $$G(\Delta,x,y,z)=\delta_1
y^{d-2}+\delta_{2}y^{d-3}z+\delta_{3}x^{d-2}+\delta_{4}x^{d-3}z
+\delta_{5}x^{d-3}y+\delta_6 xy^{d-3}.$$ If $d=3$ we take
$\Delta=(\delta_1,\delta_2,\delta_3)$ and $G= \delta_1 y+\delta_2
z+\delta_3 x$. Observe also that in \cite{PSS}, Lemma 1, it is
proved that
  for $\epsilon$-monomial curves, and hence for $d=3$, $\overline{\cP}(t)$ is always a parametrization.
  Then we consider
$\overline{H}^{**}(\Delta,t,x,y,z)=\overline{H}^{*}_{d-2}(t,x,y,z)+G(\Delta,x,y,z)$; that is
 $$\overline{H}^{**}(\Delta,t,x,y,z)=H_1(x,y,z)+t H_2(x,y,z)+G(\Delta,x,y,z).$$
Note that we are perturbing $H_1$ and hence $H_2$ keeps the required conditions on the point at infinity of $\cc$.

\parasmall

In this situation, repeating the above process with  $\overline{H}^{**}$ and $F$, instead of with $\overline{H}^{*}_{d-2}$ and $F$, we introduce
$\cSSsimplex,\cSSsimpley, \cBBsimplex,\cBBsimpley$ $\cRRsimplex,\cRRsimpley$ and  $\cPP$. So
\[ \cSSx=\resultant_y(\overline{H}^{**}(\Delta,t,x,y,1),f), \, \cSSy=\resultant_x(\overline{H}^{**}(\Delta,t,x,y,1),f),\]
and $\cBBx, \cRRx$ are the quotient and the remainder of the division of $\cSSsimplex$ by $A_1(x)$, respectively. Similarly,  for $\cBBy, \cRRy$ using $\cSSsimpley$ and $A_2(y)$. Finally, the  components of $\cPP$ are the roots of $\cBBx$ and $\cBBy$ as univariate polynomials over ${\Bbb C}[\Delta,t]$.

\parasmall

We start with some lemmas.

\parasmall

\begin{lemma}\label{lema-param-2}
The leading coefficient w.r.t. $x$ of
$\cBBx$ and the leading coefficient of  $\cBBy$ w.r.t. $y$, as polynomials in ${\Bbb C}(\Delta)[t]$, are the same up to multiplication by non-zero constants in $\Bbb C$. Furthermore, the roots are
\[  \left\{-\frac{H_1(a,b,0)+G(\Delta,a,b,0)}{H_2(a,b,0)}\right\}_{(a:b:0)\in \cc^h}. \]
\end{lemma}

\parasmall

\noindent {\bf Proof.} Let $\cBBx=q_1(\Delta,t)x- p_1(\Delta,t)$, and $\cBBy=q_2(\Delta,t)y- p_2(\Delta,t)$.
By hypothesis   $F(1,0,0)\neq 0, F(0,1,0)\neq 0$. So,  the leading coefficient of  $F$ w.r.t. $y$ is a non-zero constant; similarly   w.r.t. $x$. Thus, by well known properties on resultants (see, e.g. Lemma 4.3.1. in \cite{Winkler}), it holds that up to multiplication by a non-zero element in $\Bbb C$:
\[ \resultant_y(\overline{H}^{**}(\Delta,t,x,y,0), F(x,y,0))=(S_{1}^{\Delta})^H(\Delta,x,0,t), \]
\[ \resultant_x(\overline{H}^{**}(\Delta,t,x,y,0), F(x,y,0))=(S_{2}^{\Delta})^H(\Delta,y,0,t), \]
where $(S_{i}^{\Delta})^H$ denotes the homogenization of $S_{i}^{\Delta}$ as polynomials in ${\Bbb C}(\Delta,t)[x,y]$. Now, observe that
\[ (S_{1}^{\Delta})^H(\Delta,x,0,t)=q_1(\Delta,t) x^{d(d-2)}, (S_{2}^{\Delta})^H(\Delta,y,0,t)=q_2(\Delta,t) y^{d(d-2)}.\]
Moreover, let $F(x,y,0)$ factor as
\[ F(x,y,0)=\prod_{i=1}^{d} (\beta_i x- \alpha_i y). \]
Since $F(0,1,0)\neq 0$ then $\alpha_i\neq 0$ for all $i$. Hence, up to multiplication by non-zero constants
\[ \resultant_y(\overline{H}^{**}(\Delta,t,x,y,0), F(x,y,0))=\prod_{i=1}^{d} \resultant_y(\overline{H}^{**}(\Delta,t,x,y,0), \beta_i x- \alpha_i y)=\]
 \[ =(-1)^{d(d-2)} x^{d(d-2)} \prod_{i=1}^{d} \overline{H}^{**}(\Delta,t,\alpha_i, \beta_i, 0). \]
 Analogously,
 \[ \resultant_x(\overline{H}^{**}(\Delta,t,x,y,0), F(x,y,0))=(-1)^{d(d-2)} y^{d(d-2)} \prod_{i=1}^{d} \overline{H}^{**}(\Delta,t,\alpha_i, \beta_i, 0). \]
So, up to multiplication by non-zero constants

\vspace{1 mm}

\noindent $q_1(\Delta,t)= q_2(\Delta,t)=\prod_{i=1}^{d} \overline{H}^{**}(\Delta,t,\alpha_i, \beta_i, 0)= \prod_{i=1}^{d} (H_1(\alpha_i,\beta_i,0)+G(\Delta,\alpha_i,\beta_i,0)+t H_2(\alpha_i,\beta_i,0)).
\qed $

\parasmall

\begin{lemma}\label{lemma-param-2a} For  all $\Delta_0\in {\Bbb C}^6$,
$\degree_t(\cBBsimplex(\Delta_0,x,t))=d$ and $\degree_t(\cBBsimpley(\Delta_0,y,t))=d$.
\end{lemma}

\parasmall

\noindent {\bf Proof.} First note that $\deg_t(\cBBsimplex)\leq d$ and $\deg_t(\cBBsimpley)\leq d$. The equality follows from the last equality in the proof of Lemma \ref{lema-param-2}, and using that $H_2(\alpha_i,\beta_i,0)\neq 0$ for all $i$. \qed

\parasmall

\begin{lemma}\label{lema-param-3} There exists a non-empty  Zariski open subset $\Omega$ of ${\Bbb C}^6$ such that if $\Delta_0\in \Omega$ then $\cBBsimplex(\Delta_0,x,t)$ and $\cBBsimpley(\Delta_0,y,t)$ are primitive  w.r.t. $x$ and $y$, respectively.
\end{lemma}

\parasmall

\noindent {\bf Proof.}  We assume that $d>3$; if $d=3$ the reasoning is analogous.
 Let us assume that
\[ \cBBx=D(\Delta,t) L(\Delta,x,t), \]
with $\deg_t(D)>0$.
Then,
\[ \cSSx=D(\Delta,t) L(\Delta,x,t)A_1(x)+\cRRx. \]
By Lemma \ref{lema-param-2}, we know how the roots of
$D(\Delta,t)\in {\Bbb C}(\Delta)[t]$ are. Now for each root $t_0$
of $D$ (say that $t_0$ is defined by $P:=(a:b:0)\in\cc^h$),
$\deg_x(\cSSsimplex(\Delta,x,t_0))=\deg_x(\cRRsimplex(\Delta,x,t_0))\leq
d(d-2)-2$. Let  ${\cal D}(t_0)$ be the projective curve defined by
$\overline{H}^{**}(\Delta,t_0,x,y,z)$ over the algebraic closure
$\Bbb F$ of ${\Bbb C}(\Delta)$. Then,  ${\cal D}(t_0)$ and $\cc^h$
intersect at infinity  at an additional point different from $P$,
or the multiplicity of intersection of both curves at $P$ is at
least two. We analyze each case. But first we introduce some
additional notation. We express $F,H_1,$ and $H_2$ as
\[ F(x,y,z)=f_0(x,y)+f_1(x,y)z+\cdots+f_{d}(x,y) z^d,\] \[ H_i(x,y,z)=h_{i,0}(x,y)+h_{i,1}(x,y)z+\cdots+h_{i,d-2}(x,y)z^{d-2}, \]
where $f_j,h_{i,j}$ are homogeneous of degree $d-j$ and $(d-2)-j$, respectively. Moreover, we denote by $F^x, H_{i}^{x},  \overline{H}^{**,x},
f_{j}^{x},h_{i,j}^{x}$ the corresponding partial derivative w.r.t. $x$; similarly w.r.t. $y$ and $z$.

Let us assume that $Q\in {\cal D}(t_0)$, with $Q=(n:m:0)\neq P$.
This is equivalent  to
\[ \delta_1 C_1+\delta_3 C_3
+\delta_5 C_5+\delta_6 C_6
=C_0, \]
where
\[ \begin{array}{ll} C_1=b^{d-2}h_{2,0}(Q)-m^{d-2}h_{2,0}(P), &
C_3=a^{d-2}h_{2,0}(Q)-n^{d-2}h_{2,0}(P), \\
C_5= a^{d-3}bh_{2,0}(Q)-n^{d-3}m h_{2,0}(P), &
C_6=a b^{d-3}h_{2,0}(Q)-n m^{d-3} h_{2,0}(P), \\
C_0=h_{1,0}(P)h_{2,0}(Q)-h_{2,0}(P)h_{1,0}(Q). &
\end{array} \]
 Observe that $h_{2,0}(Q)\neq 0$, $h_{2,0}(P)\neq 0$.  Let us see that all $C_i, i>0,$ can not vanish simultaneously. Let $C_1=C_3=C_5=C_6=0$.
We assume that $a\neq 0$. If $a=0$ then $b\neq 0$, and the reasoning is similar. From $C_3=0$ one has that $n\neq 0$. So
\[ P=(a:b:0)=(a^{d-2}:ba^{d-3}:0)=(a^{d-2}h_{2,0}(Q):ba^{d-3}h_{2,0}(Q):0)=\]\[=(n^{d-2}h_{2,0}(P):n^{d-3}m h_{2,0}(P):0)=(n:m:0)=Q, \]
 which is a contradiction. Therefore,  if $V_1$ is the hyperplane in ${\Bbb C}^6$ defined by $\delta_1 C_1+\delta_3 C_3
+\delta_5 C_5+\delta_6 C_6=C_0$, for all $\Delta_0$ in ${\Bbb C}^6 \setminus V_1$ this case does not happen.

Let us assume that the multiplicity of intersection of ${\cal
D}(t_0)$ and $\cc^h$ at $P$ is at least two. Since $\cc^h$ does
not have singularities at infinity, this implies that both curves
have the same tangent at $P$. This is equivalent to demand
\[ (F^x(P):F^y(P):F^z(P))=(\overline{H}^{**,x}(\Delta,t_0,P): \overline{H}^{**,y}(\Delta,t_0,P): \overline{H}^{**,z}(\Delta,t_0,P)). \]
By hypothesis $a b\neq 0$. So, by Euler's formula and taking into
account that $P$ is at infinity, the condition is equivalent to
$$F^x(P)\overline{H}^{**,z}(\Delta,t_0,P)=F^z(P)\overline{H}^{**,x}(\Delta,t_0,P).$$
That is equivalent to
\[ \delta_1 C_1+\delta_2 C_2+\delta_3 C_3+
\delta_4 C_4+\delta_5 C_5+\delta_6 C_6
=C_0, \]
where
\[ \begin{array}{l} C_1=b^{d-2}(f_1(P)h_{2,0}^{x}(P)-f_{0}^{x}(P)h_{2,1}(P)), \\
C_2=b^{d-3}f_{0}^{x}(P) h_{2,0}(P), \\
C_3=a^{d-3}(af_1(P)h_{2,0}^{x}(P)-af_{0}^{x}(P)h_{2,1}(P)-(d-2)f_{1}(P)h_{2,0}(P)), \\
C_4=a^{d-3}f_{0}^{x}(P) h_{2,0}(P), \\
C_5=a^{d-4}b(-f_{0}^{x}(P)h_{2,1}(P)a +f_{1}(P)h_{2,0}^{x}(P)a   -(d-3)f_1(P)h_{2,0}(P)),  \\
C_6= b^{d-3}(f_1(P)h_{2,0}^{x}(P)a-f_1(P)h_{2,0}(P)-f_{0}^{x}(P)h_{2,1}(P)a), \\
C_0=f_{0}^{x}(P)(h_{2,1}(P)h_{1,0}(P)-h_{1,1}(P)h_{2,0}(P))+f_1(P)(h_{1,0}^{x}(P)h_{2,0}(P)-h_{2,0}^{x}(P)h_{1,0}(P)).
\end{array} \]
Let us see that $C_i, i>0$, cannot vanish simultaneously. Let $C_1=\cdots=C_6=0$. Since $h_{2,0}(P)\neq 0$, and $ab\neq 0$, one has that
\[ \begin{array}{ll} C_1=0 \Rightarrow & f_1(P)h_{2,0}^{x}(P)=f_{0}^{x}(P)h_{2,1}(P), \\
C_2 =0 \Rightarrow &f_{0}^{x}(P)=0, \\
C_6=0 \Rightarrow & f_1(P)(h_{2,0}^{x}(P)a-h_{2,0}(P))=0.
\end{array} \]
Note that $f_{1}(P)\neq 0$, since otherwise it would imply that $F^z(P)=0$ and using that $F^x(P)=f_{0}^{x}(P)=0$ and that (by Euler's formula) $F^{y}(P)=0$, one would deduce that $P$ is a singularity of $\cc^h$ which is excluded by hypothesis. So, the first and second equalities imply that $h_{2,0}^{x}(P)=0$ and this yields to (using the last equality) $h_{2,0}(P)=0$ which is a contradiction. Therefore, out of the hyperplane defined in ${\Bbb C}^6$ by $\sum_{i=1}^{6} C_i \delta_{i}=C_0$, this case cannot happen.

For each point of $\cc$ at infinity we generate the hyperplanes
described above and corresponding to each one of the two cases.
Let $V$ be the union of all of them, and let  $\Omega_1={\Bbb
C}^{6} \setminus V$. Repeating the same reasoning with
$\cBBsimpley$ (note that $G$ is symmetric in terms of $x$ and
$y$), we get $\Omega_2$. Finally, let $\Omega=\Omega_1\cap
\Omega_2$. \qed

\parasmall

Now, the next theorem follows directly.

\parasmall

\begin{theorem}\label{theorem-param} There exists a non-empty  Zariski open subset $\Omega$ of ${\Bbb C}^6$ such that if $\Delta_0\in \Omega$ then $\cPPsimple(\Delta_0,t)$ is  a rational parametrization of a rational curve of degree at most $d$.
\end{theorem}

\parasmall

\noindent {\bf Proof.} Taking $\Omega$ as in Lemma \ref{lema-param-3}, we ensure that $\cPPsimple(\Delta_0,t)$ is a rational parametrization. By Lemmas \ref{lema-param-2} and \ref{lemma-param-2a} we get that the degree of the curve is at most $d$. \qed

\parasmall

\begin{remark}
    Let $\overline{\cal H}^{**}_{\Delta_0}$ be the linear system of $(d-2)$-degree curves defined by $\overline{H}^{**}(\Delta_0,t,x,y,z)$.  If no perturbation is needed, i.e. $\Delta=\vec{0}$, then  $\overline{\cal H}^{**}_{\vec{0}}=\overline{\cal H}^{*}_{d-2}$, and hence it is generated by the effective (exact) divisor $\sum_{i=1}^{s} r_i Q_i+\sum_{i=1}^{d-3} P_i$.
    Now, if we identify (as usual)   ${\Bbb C}^2$ with ${\Bbb R}^4$ and we consider the perturbing parameters $\delta_i$ as real variables, it holds that for each $Q_i$ (similarly of $P_i$) there exists $\rho(Q_i)>0$ such that for almost every element $(Q_{i}^{*},\Delta_0) $ in the open Euclidean  disk of ${\Bbb R}^{10}$, of center $(Q_i,\vec{0})$ and radius $\rho(Q_i)$,
        $\overline{\cP}(\Delta_0,t)$ is a parametrization and $Q_{i}^{*}$ is an $\epsilon$-point of $\epsilon$-multiplicity (at least) $r_i$  of
the curves defined by $H_1(x,y,z)+G(\Delta_0,x,y,z)$ and $H_2(x,y,z)$ (recall that $\overline{H}^{**}=H_1+G+t H_2$); i.e. of the generating curves of  $\overline{\cal H}^{**}_{\Delta}$.   This can be seen by applying Theorem \ref{theorem-param}, taking $\delta_i$ small enough to ensure that $\|H_{1}\|=\|H_1+G(\Delta_0,x,y,z)\|$, and noting that if $M$ is any of the derivatives of $H_1+G$ and $H_2$ involved in the $\epsilon$-multiplicity,  then   $\| M \|_{2}^{2}$ is a  continuous function that  vanishes at $(Q_i,\vec{0})$. \hfill \qed
\end{remark}

\parasmall

\parasmall

Finally, and before outlining the algorithm, we briefly describe
how to proceed with the selection and computation of the (affine
simple) $\epsilon$-points $P_i$. We first observe that, in
general, an $\epsilon$-point can be computed by solving
$\{f(x,y)=0,\alpha x+\beta y=\rho\}$, where $\alpha, \beta,\rho\in
\Bbb C$,  under fixed precision $\epsilon \|f\|$. However, we are
intersected in working with either real $\epsilon$-points  or
pairs of conjugate complex points. We can always compute all
points, but at most one, in pairs of conjugate complex points. For
choosing real points one can always analyze the roots of the
discriminant of $f$ (see Theorem 7.7 in \cite{SWP}). On the other
hand we have observed, in our examples, that taking (when
possible) the simple $\epsilon$-points as (affine)
$\epsilon$-ramification points (see Def. \ref{def--mult-e-point})
the error distance between the original curve and the output curve
decreases. So we tend to use first such points. Finally, one has
to take care of the fact that a chosen $\epsilon$-point can be too
close (i.e. in the same cluster) to an $\epsilon$-singularity or
to a previously computed $\epsilon$-point, and hence identifiable
with it. To avoid this, whenever a  new simple  $\epsilon$-point
is computed we check whether it belongs to the cluster of the
others points.

\parasmall

 The above process provides the following approximate
parametrization algorithm for deciding whether  a real $\epsilon$-irreducible (with proper degree) plane algebraic
curve $\cal C$ is $\epsilon$--rational, and in the affirmative case,
compute an approximate parametrization. Recall that we assume that $\cc$ does not have exact singularities at infinity, and that
$(0:1:0),(1:0:0)\notin \cc^h$. If this last condition fails, one may consider an affine orthogonal change of coordinates to achieve the requirement.

\parasmall

\noindent
\noindent \underline{\sf Approximate Parametrization Algorithm}
\begin{itemize}
\item \underline{\sf Given}  a tolerance $\epsilon >0$ and an
$\epsilon$--irreducible polynomial    $f(x,y)  \in {\Bbb R}[x,y],$
of proper degree $d>2$ (for $d=1$ it is trivial, if $d=2$ apply
\cite{PSS}), without exact singularities at infinity, not passing
through $(0:1:0), (1:0:0)$, and defining a real plane algebraic
curve ${\cal C}$; let $F(x,y,z)$ be the homogenization of $f$.
\item \underline{\sf Decide} whether $\cal C$ is
$\epsilon$-rational and in the affirmative case \item
\underline{\sf Compute} a rational parametrization $\overline{\cal
P}(t)$ of a curve $\overline{\cal C}$ close to ${\cal C}$.
\end{itemize}
\begin{itemize}
\item[(1)] {\sf Compute} the cluster decomposition
$\{\cluster_{r_i}(Q_i)\}_{i=1,\ldots,s}$ of $\sing(\cc)$; say
$Q_i=(q_{i,1}:q_{i,2}:1)$. \item[(2)] If
$\sum_{i=1}^{s}r_i(r_i-1)\neq (d-1)(d-2),$  {\sf RETURN} ``{\it
$\cal C$ is not (affine) $\epsilon$-rational}". If $s=1$ one may
apply the algorithm in \cite{PSS}. \item[(3)] {\sf Determine} the
linear system  $\overline{\cal H}_{d-2}$  of degree $(d-2)$ given
by the divisor $\sum_{i=1}^{s} r_i Q_i$. \item[(4)] {\sf Compute}
$(d-3)$  $\epsilon$--ramification points $\{P_j\}_{1\leq j\leq
d-3}$ of ${\cal C}$; if there are not enough
$\epsilon$-ramification points, complete with simple
$\epsilon$-point. Take the points over $\Bbb R$, or as conjugate
complex points. After each point computation check that it is not
in the cluster of the others (including the clusters of $Q_i$); if
this fails take a new one. Say $P_i=(p_{i,1}:p_{i,2}:1)$.
\item[(5)] {\sf Determine} the linear subsystem $\overline{\cal
H}^{\,*}_{d-2}$ of $\overline{\cal H}_{d-2}$ given by the divisor
$\sum_{i=1}^{d-3}P_i$. Let $H^*(t,x,y,z)=H_1(x,y,z)+tH_2(x,y,z)$
be its defining polynomial.
 \item[(6)]If[$\gcd(F(x,y,0),H_1(x,y,0))\neq 1$] and
[$\gcd(F(x,y,0),H_2(x,y,0))\neq 1$]  replace $H_2$ by
$H_2+\rho_1x^{d-2}+\rho_2y^{d-2}$, where $\rho_1, \rho_2$ are real
and strictly smaller than $\epsilon$.  Say that
$\gcd(F(x,y,0),H_2(x,y,0))=1$; similarly in the other case.
\item[(7)] Set $\delta_1=\cdots=\delta_6=0$. \item[(8)] If $d>3$
then $g:=\delta_1
y^{d-2}+\delta_{2}y^{d-3}+\delta_{3}x^{d-2}+\delta_{4}x^{d-3}+\delta_{5}x^{d-3}y+\delta_6
xy^{d-3}$ else $g:=\delta_1 y+ \delta_2 + \delta_3 x$. \item[(9)]
$ {S}_1(x, t)=\resultant_y(H^{*}(x,y,1)+g, {f})$ and   $ {S}_2(y,
t)=\resultant_x(H^{*}(x,y,1)+g, {f}).$
 \item[(10)] $A_1(x)=\prod_{i=1}^{s}(x-{q}_{i,1})^{r_i(r_i-1)}
\prod_{i=1}^{d-3}(x-{p}_{i,1}), \\
A_2(y)=\prod_{i=1}^{s}(y-{q}_{i,2})^{r_i(r_i-1)}\prod_{i=1}^{d-3}(y-{p}_{i,2})$.
\item[(11)]  For $i=1,2$ {\sf compute}  the   quotient $B_i$  of $
{S}_i$ by $A_i$ w.r.t.  either $x$ or $y$. \item[(12)] If the
content of  $B_1$ w.r.t $x$  or the content of $B_2$ w.r.t. $y$
does  depend on $t$, take $\{\delta_1,\ldots,\delta_6\}$ as  small
real numbers (strictly smaller than $\epsilon$) and go to Step 8.
\item[(13)]  {\sf Determine} the root $\overline{p}_1(t)$ of $B_1$
as a polynomial in $x$ and  the root $\overline{p}_2(t)$ of $B_2$
as a polynomial in $y$. \item[(14)] {\sf RETURN}  $\overline{\cal
P}(t)=(\overline{p}_1(t),\overline{p}_2(t))$.
\end{itemize}

The next theorem states the main properties of the curve output by the algorithm. But first, we need the following
technical lemma.

\parasmall

\begin{lemma}\label{lemma-para-ptos-inf}
Let $\Bbb L$  be the algebraic closure of ${\Bbb C}(t)$, and $\cc_1, \cc_2$ two plane projective curves over $\Bbb L$ with defining polynomials $G_1(x,y,z),G_2(x,y,z)\in {\Bbb C}[t][x,y,z]$, respectively. If there exist $K,W,L\in {\Bbb C}[t][x,y,z] $ such that $KG_1+W G_2=z L$, and
\begin{itemize}
\item[(1)] $G_1(x,y,0)G_2(x,y,0)\neq 0$,
\item[(2)] $\gcd(G_1(x,y,0),G_2(x,y,0))=1$,
\end{itemize}
then either $z$ divides $K$ and $W$ or there exist $U_1,U_2,U_3\in {\Bbb C}[t][x,y,z]$ such that
\[ L=U_1 G_1(x,y,0)+U_2 G_2(x,y,0) +z U_3. \]
\end{lemma}

\parasmall

\noindent {\bf Proof.} If $z$ divides $K$, then $z$ divides $WG_2$, and by (2) $z$ divides $W$. So let us assume that $z$ does not divides $K$, and let us denote by $G_{i}^{0}$ the polynomial $G_i(x,y,0)$; similarly with $K^0,W^0$. Then, $K^0G_{1}^{0}+W^0G_{2}^{0}=0$. Since $G_{i}^{0}\neq 0$ and $\gcd(G_{1}^{0},G_{2}^{0})=1$, then
$G_{1}^{0}$ divides $W^0$ and $G_{2}^{0}$ divides $K^0$. Let $K^0=\Delta_1 G_{2}^{0}, W^0=\Delta_2 G_{1}^{0}$. So $(\Delta_1+\Delta_2) G_{1}^{0}G_{2}^{0}=0$, and since $G_{i}^{0}\neq 0$, one gets $\Delta_1+\Delta_2=0$. Now,  we write
\[ K=K^0+z\overline{K}, W=W^0+z\overline{W}, G_{i}=G_{i}^0+z\overline{G}_{i},
 \]
where $\overline{K}, \overline{W},\overline{G}_i \in {\Bbb C}[t][x,y,z]$. Then,
$KW_1+WG_2=z(G_{1}^{0}\overline{K}+G_{2}^{0}\overline{W}+z(\overline{K}\, \overline{G}_1+\overline{W}\, \overline{G}_2))$.
 \qed

\parasmall

\begin{theorem} \label{T-infinity}
The rational curve $\overline{\cc}$, output by the algorithm, and $\cc$
have the same points at infinity, and $\deg(\overline{\cc})\leq \deg(\cc)$.
\end{theorem}

\parasmall

\noindent {\bf Proof.} The fact on  the degree follows from Theorem \ref{theorem-param}.
 For the reasoning of the rest of the proof, we can assume w.l.o.g. that
no perturbation $\Delta$ is required (i.e. $\Delta=\vec{0}$) in the execution of the algorithm.
Let $\overline{H}^{*}_{d-2}(t,x,y,z), F(x,y,z), S_1(x,t), S_2(y,t), A_1(x), A_2(y), B_1(x,t), B_2(y,t)$, $R_1(x,t),$ and  $R_2(y,t)$ be defined as above. Let
$B_1:={q}_1(t) x-{p}_1(t)$,  $B_2:={q}_2(t) y-{p}_2(t)$, and recall that $R_i$ is the remainder of the division of $S_i $ by $A_i$. Furthermore, by Lemma \ref{lema-param-2}, ${q}_1(t)=\lambda {q}_2(t)$, with $\lambda\in {\Bbb C}^*$.
By Lemma \ref{lemma-param-2a}, $\deg_t(B_1)=\deg_t(B_2)=\deg(F)=d$, and, by Lemma \ref{lema-param-3}, $\gcd({q_1},{p}_1)= \gcd({q}_2,{p}_2)=1$. So,
\[ \overline{\cP}^{H}(t):=(\lambda^{-1} p_1(t):p_2(t):q_2(t)) \]
parametrizes the projective closure of $\overline{\cc}$. Furthermore, since $\deg(p_i)\leq \deg(q_2)$, then all points of $\overline{\cc}$ at infinity are reachable by $\overline{\cP}^{H}(t)$ (see \cite{Sendra}).
In addition, we note that
\[  \deg_{\{x,y,z\}}(\overline{H}_{d-2}^{\,*})=d-2,\,
\deg(A_j)=d(d-2)-1,\,
\deg_{\{x,y\}}(R_j)\leq d(d-2)-2.\] Moreover, if
$m^{H}(x,y,z,w)$ denotes the homogenization of $m(x,y,w)$  as a
polynomial in ${\Bbb C}[w][x,y]$, we have that
$$\overline{S}^H_1(x, z, t)=\resultant_y(\overline{H}_{d-2}^{*}(t,x,y,z),F(x,y,z))=B^H_1(x,z,
t)A^H_1(x,z)+R^H_1(x,z,t)z^{n_1},\,$$
$$\overline{S}^H_2(y, z, t)=\resultant_x(\overline{H}_{d-2}^{*}(t,x,y,z),F(x,y,z))= {B}^H_2(y,z,
t)A^H_2(y,z)+R^H_2(y,z,t)z^{n_2},\,$$
where $n_j+\deg(R_j^H)= d(d-2),\,j=1,2$. So $n_j\geq 2$. Also, we denote by $\overline{\cc}_\infty$ and $\cc_\infty$ the set of points at infinity $\overline{\cc}$ and $\cc$ respectively.
By resultant properties, there exist polynomials $M_1,N_1,M_2,N_2\in {\Bbb C}[t,x,y,z]$ such that
\[ M_i \overline{H}_{d-2}^{*} +N_i F= \overline{S}^H_i, \,\, i=1,2. \]
So,
\[ y A_{2}^{H} S_{1}^H-\lambda x A_{1}^{H} S_{2}^{H}= z A_{1}^{H}A_{2¨}^{H}(\lambda x p_2- y p_1)+ z^{n_3} R_{3}, \]
where $n_3\geq 2$ and $R_3$ a polynomial; namely $z^{n_3}R_3= yA_2 z^{n_1}R_{1}^{H}-\lambda x A_1 z^{n_2}R_{2}^{H}$.
On the other hand,  if $K=yA_{2}^{H}M_1-\lambda x A_{1}^{H} M_2$ and, $W=yA_{2}^{H} N_1-\lambda x A_{1}^{H} N_2$, then
\[ y A_{2}^{H} S_{1}^H-\lambda x A_{1}^{H} S_{2}^{H}=K(x,y,z,t) (H_1 +t H_2)+ W(x,y,z,t) F. \]
Therefore, $z$ divides the right hand side of the above equation. We now check that $H_1+tH_2$ and $F$ satisfy the hypothesis of Lemma  \ref{lemma-para-ptos-inf}. Since $F$ is irreducible and non-linear, $F(x,y,0)\neq 0$. Moreover, if $H_1(x,y,0)+tH_2(x,y,0)=0$ then $H_2(x,y,0)=0$ and this implies that ${\cal D}_2$ contains all the points at infinity of $\cc^h$, which is a contradiction. Finally, if $\gcd(H_1(x,y,0)+t H_2(x,y,0),F(x,y,0))\neq 1$, then $\gcd(H_2(x,y,0),F(x,y,0))\neq 1$, and this  implies that ${\cal D}_2$ and $\cc^h$ share points at infinity.  Therefore, applying Lemma \ref{lemma-para-ptos-inf},
one deduces that either
 there exist polynomials $M_3,N_3\in {\Bbb C}[t][x,y,z]$ such that
\[ M_3 \overline{H}^{*}_{d-2}+N_3 F= A_{1}^{H}A_{2¨}^{H}(\lambda x p_2- y p_1)+ z^{n_4} R_{3}, \]
where $n_4\geq 1$, or
 there exist polynomials $U_1,U_2,U_3\in {\Bbb C}[t][x,y,z]$  such that
 \[ U_1  \overline{H}^{*}_{d-2}(t,x,y,0)+U_2  F(x,y,0)+z U_3 =A_{1}^{H}A_{2¨}^{H}(\lambda x p_2- y p_1)+ z^{n_4} R_{3}.\]
In this situation, using $\overline{\cc}_\infty\subset
\overline{\cP}^{H}({\Bbb C})$, we first observe that
$\card(\overline{\cc}_\infty)$ is less or equal to the number of
different roots of $q_2(t)$ and, by Lemma \ref{lema-param-2}, this
number is less or equal to $\card({\cc}_\infty)$. So,
$\card(\overline{\cc}_\infty)\leq \card({\cc}_\infty)$.  Now, we
prove that $\cc_\infty\subset \overline{\cc}_\infty$, from where
one concludes the proof. Let $P=(x_0:y_0:0)\in \cc_\infty$, and
let $t_0$ be the root of $q_2$ generated by $P$ (see Lemma
\ref{lema-param-2}). So,
$\overline{H}_{d-2}^{*}(t_0,x_0,y_0,0)=F(x_0,y_0,0)=0$. Applying
the corresponding equality above, and using that $n_4\geq 1$, we
get
\[ A_{1}^{H}(x_0,0)A_{2}^{H}(y_0,0) (\lambda x_0 p_2(t_0)-y_0 p_1(t_0))=0. \]
Moreover, since $(1:0:0),(0:1:0)\not\in \cc^h$ then $x_0y_0\neq 0$, and hence $ A_{1}^{H}(x_0,0)A_{2}^{H}(y_0,0)\neq 0$. So,
$\lambda x_0 p_2(t_0)=y_0 p_1(t_0)$. In addition, $p_1(t_0)p_2(t_0)\neq 0$ because $\gcd(q_2,p_1)=1=\gcd(q_2,p_2)$. Therefore,
\[ \overline{\cP}^{H}(t_0)=(\lambda^{-1}p_1(t_0):p_2(t_0):0)=\]
\[=(y_0\lambda^{-1}p_1(t_0):y_0p_2(t_0):0)=(x_0p_2(t_0):y_0p_2(t_0):0)=(x_0:y_0:0)=P.  \]
 \qed

\section{Displaying Examples.}\label{sec-examples}

In this section we present several examples (the degrees are 5,6,7)  to illustrate the
algorithm.
 These examples have been computed  in Maple.

\parasmall

\begin{example}\label{example-degree-5}
Let  $\epsilon= 0.01$ and ${\cal C}$ the curve of proper
degree $5$ defined by the polynomial (see Fig.\ref{curves5}):

 \vspace{1 mm}

\parbox{15cm}{\small
$f(x,y)=
0.006521014507x^4+0.006521014507x^3y^2-0.3174429862x^3+0.006521014507x^4y
+0.03536521618y^4+
0.008903520149x^2y^3-0.1541837293y^3-0.3561209555x^2-0.2351465855y^2+0.01517989182xy^4+
0.006177658243y^5+0.006521014507x^5-0.6503293396xy+0.006521014507x^3y-0.6965951291xy^2+
0.1751383118x^2y^2+0.1487535027xy^3-1. x^2y+0.0000658688334.$
}

 \parasmall

\noindent First we
compute the $\epsilon$-singularities of ${\cal C}$, obtaining
the  $\epsilon$-singular locus $\sing({\cal C})={\cal S}_{1}\cup
{\cal S}_{2}\cup {\cal S}_{3}$:
\[
\begin{array}{l}
\begin{array}{ll}
{\cal S}_1=\{\hspace*{-0.3cm}& Q_{1}=(-3.999854219, 2.000094837),
          \, Q_{2}=(0., 0.), \\ \, & Q_{3}=(0.9998153818, -2.999388343),\\
           \,& Q_{4}=(-2.001190360+0.05414244305i, 3.001898191-0.08039416354i),\\
\,&Q_{5}=(-1.980207988, 3.002780607),
 \, Q_{6}=(-2.019931003, 2.997118979),\\
 \, &Q_{7}=
(-2.001190360-0.05414244305i, 3.001898191+0.08039416354i)\},
\end{array}
\\
\begin{array}{ll}
{\cal S}_2=\{\hspace*{-0.3cm}& Q_{8}=(-2.000000001, 3.000000001)\},
\end{array}
\\
\begin{array}{ll}
{\cal S}_3=\emptyset.
\end{array}
\end{array}
\]

\parasmall

\noindent  Moreover, the cluster decomposition of the singular
locus is (see Fig. \ref{clusters5}, Left):

 \parasmall

\parbox{15.5cm}{\small  $
\cluster_2(Q_{1})=\{Q_1\},\,\cluster_2(Q_{2}) =\{Q_2\},\,\,\cluster_2(Q_{3}) = \{Q_3\}\,\,\,\mbox{and}\\
\cluster_3(Q_{8})= \{{Q_4,Q_5,Q_6,Q_7,Q_8}\}.$}

\parasmall

\noindent We observe that ${\cal C}$ is $\epsilon$-rational.
Following Step $4$ in the algorithm we obtain  two
$\epsilon$-ramification points, namely
   $P_{1}=(3.437938023,4.260660564),$
  $P_{2}=(7.712891931,1.573609575)$.
We note that these  points are not in the cluster of each other
and they are not in the clusters of  the  cluster decomposition of
the singular locus (see Fig. \ref{clusters5}, Right).
\begin{center}
\begin{figure}
\centerline{
\psfig{figure=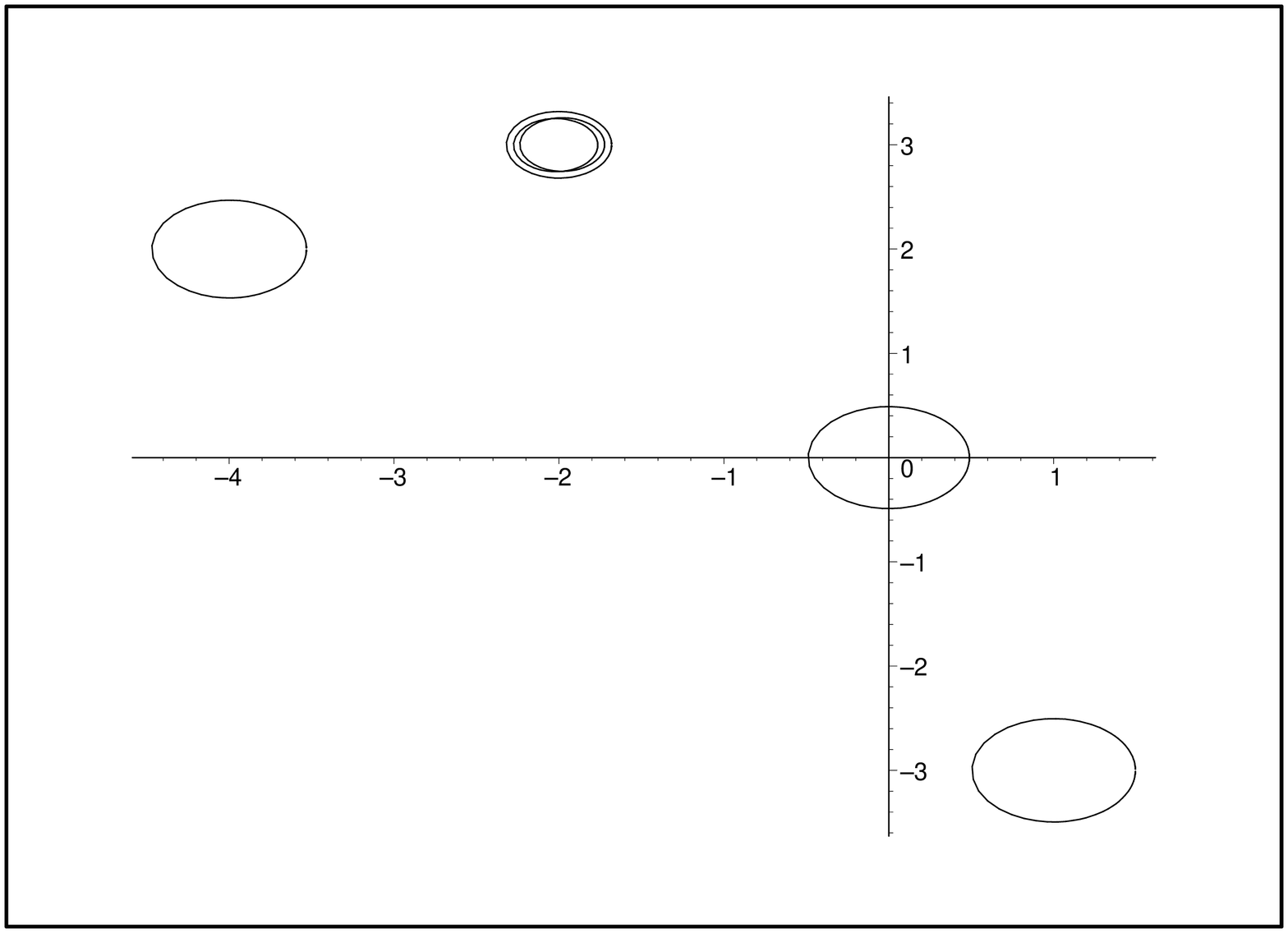,width=5.8cm,height=5.8cm,angle=270}
\hspace*{1 cm}
\psfig{figure=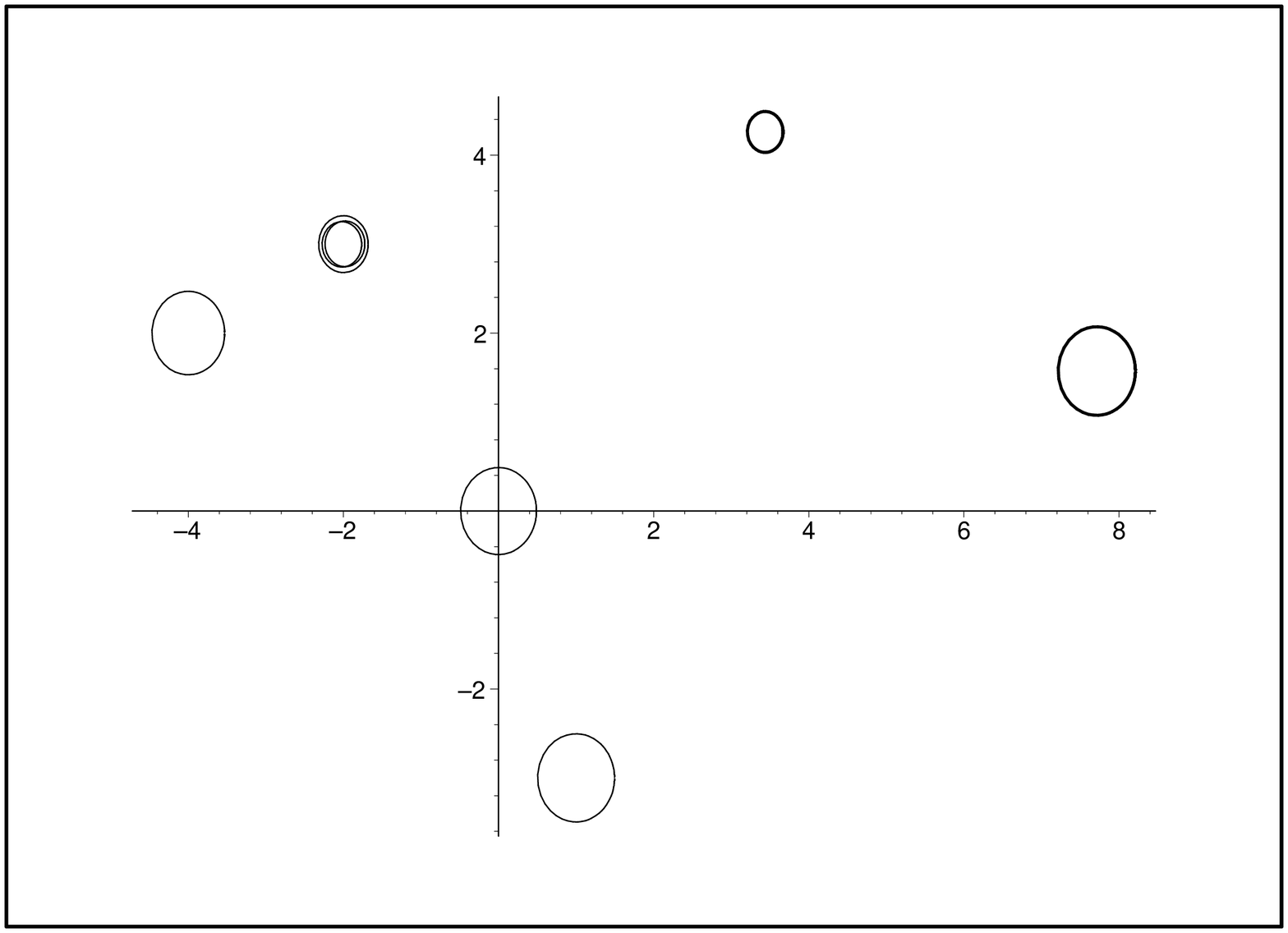,width=5.8cm,height=5.8cm,angle=270 }
 }
\caption{{\small {\sf Left:} Cluster Decomposition of the Singular
Locus.  {\sf Right:}  Cluster Decomposition of the Singular
Locus with two $\epsilon$-ramification points. }}
\label{clusters5}
\end{figure}
\end{center}
 Finally,   the algorithm   outputs the parametrization
$\overline{\mathcal
P}(t)=(\frac{\overline{p}_1(t)}{\overline{q} (t)},\frac{\overline{p}_2(t)}{\overline{q} (t)})$
 where (see Fig.
\ref{curves5} to compare the input and the output curves):

\parasmall

\noindent\parbox{14cm}{\small
 $\overline{p}_1(t)=-0.5918689071 \,10^{-29}
 (0.7256428750 \,10^{579}t+0.1009796140
 \,10^{581}t^4+0.4757134093\, 10^{580}t^5+0.3531628351 \, 10^{580}t^2+
 0.8491037424 \, 10^{580}t^3+0.5883163866\, 10^{578})
 $}

\noindent
\parbox{14cm}{\small
 $\overline{p}_2(t)= 0.3851669500\, 10^{-31}(0.1621127956\, 10^{583}t^3+0.1491645111\, 10^{582}t+
 0.6997743856\, 10^{582}t^2+0.8444468165\, 10^{582}t^5+0.1252710479\, 10^{581}+0.1858263849\, 10^{583}t^4)$}

\noindent\parbox{14cm}{\small $\overline{q}(t)=0.1265532998\,
10^{551}t^3+0.1372217100\, 10^{551}t^4+0.1260572385\,
10^{549}+0.1356321818\, 10^{550}t+ 0.5851539780\,
10^{550}t^2+0.5967959572\, 10^{550}t^5$}

\begin{center}
\begin{figure}[ht]
\centerline{ \psfig{figure=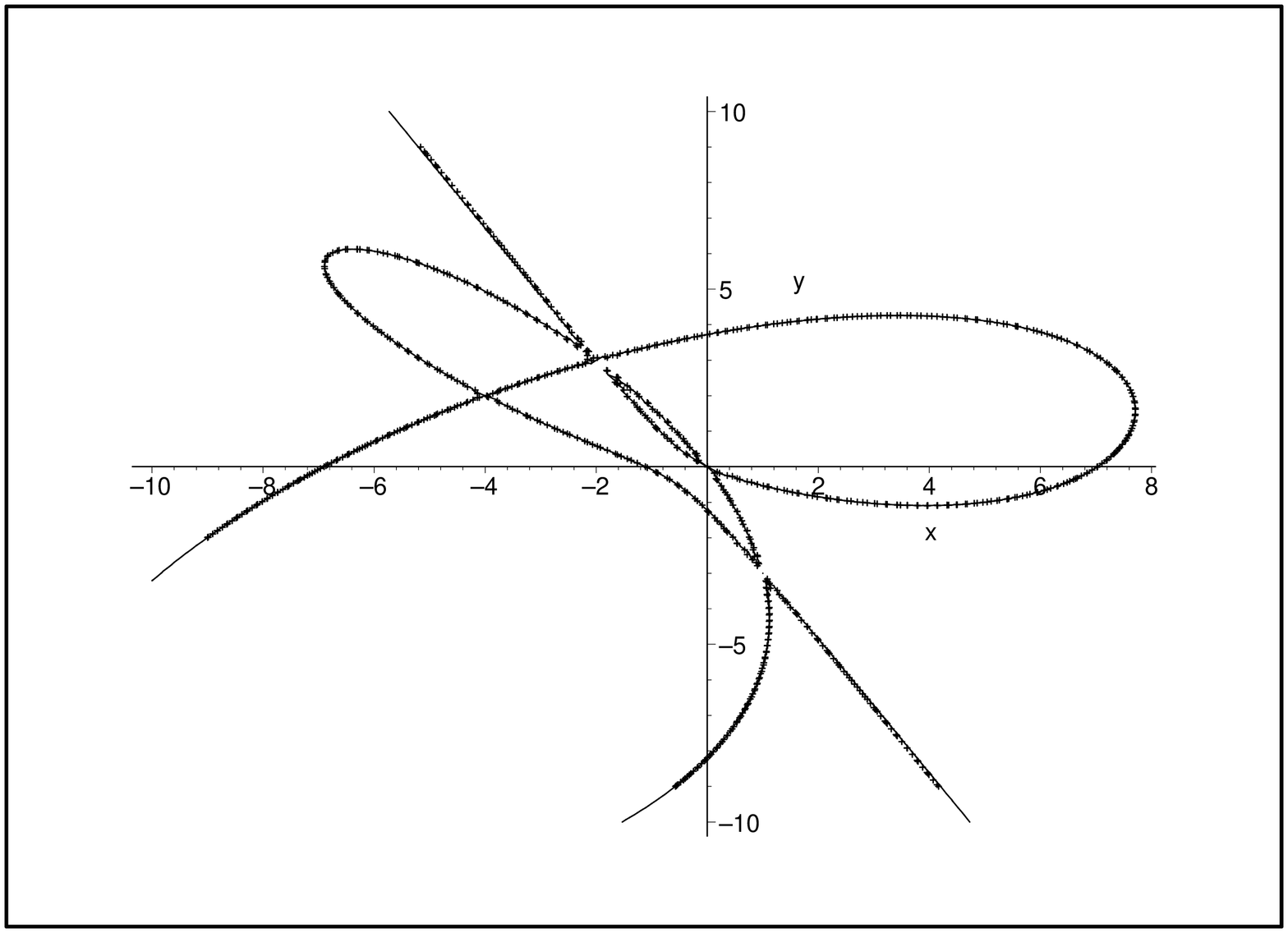,width=6.4cm,height=5.6cm,angle=270}
\psfig{figure=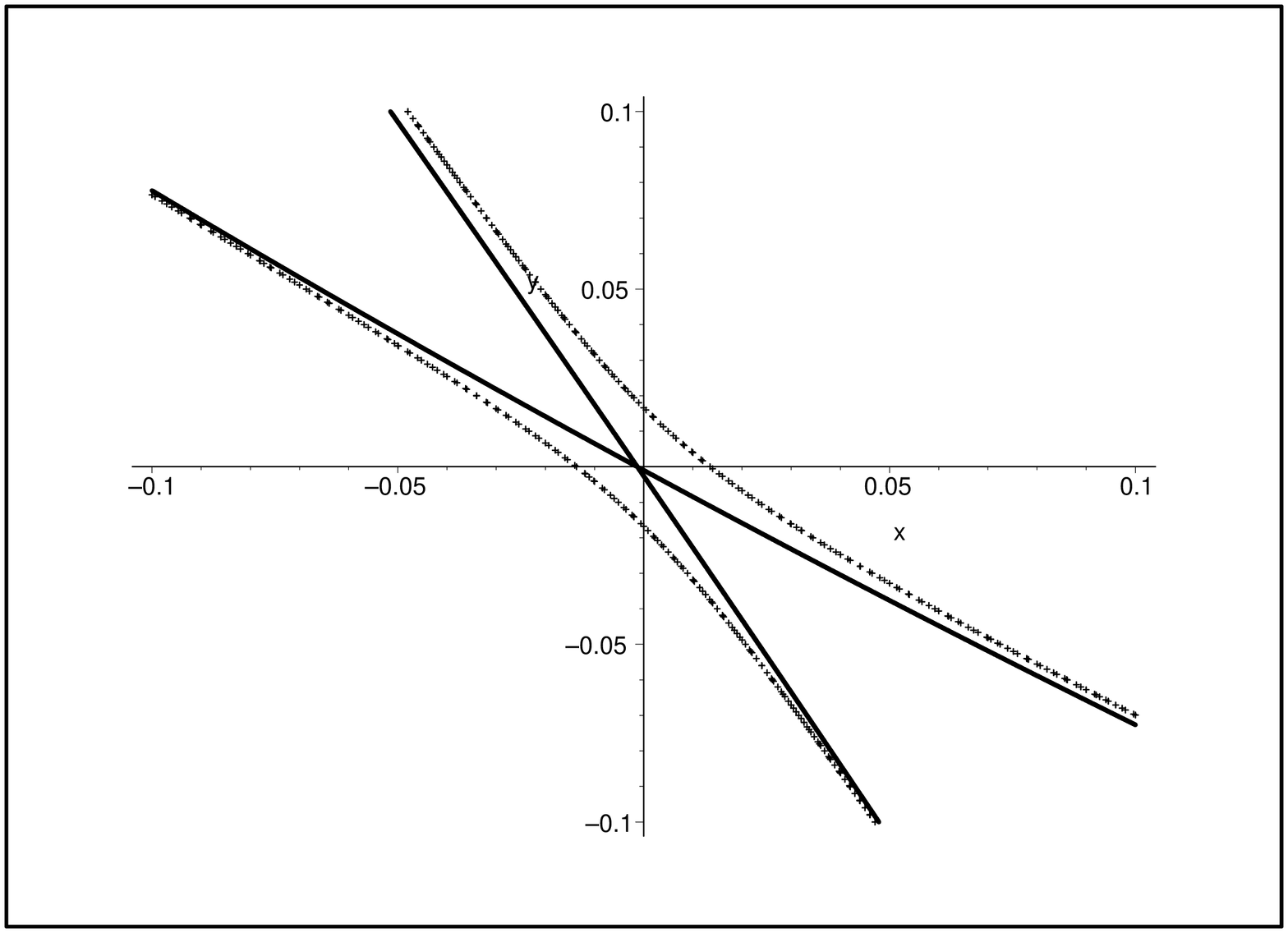,width=6.4cm,height=5.6cm,angle=270}
  } \caption{{\sf Left:} Input (in dots) and output curve in Example \ref{example-degree-5}. {\sf Right:}
A zoom at the origin
 }
\label{curves5}
\end{figure}
\end{center}

 We note that the algorithm did no require
perturbing $\overline{H}^{**}_{d-2}$.\hfill \qed

\end{example}

\parasmall

\begin{example}\label{example-degree-6}
Let $\epsilon= 0.004$ and  ${\cal C}$ the curve of proper
degree $6$ defined by the polynomial (see Fig.\ref{curves6}):
 \vspace{1 mm}

\parbox{15cm}{\small
 $f(x,y)= - 0.5499999998\,x + 0.9999999999\,y + 0.00002677376171 +
0.006666666664\,x^{2}\,y^{3}  \mbox{} + 0.006666666664\,x^{4}\,y +
0.3799999999\,x\,y^{2} - 0.4133333332\,x\,y^{3}  \mbox{} +
0.006666666664\,x^{5}\,y - 0.01999999999\,x\,y^{5} -
0.1066666667\,x^{3}\,y  \mbox{} - 0.07000000000\,x^{2}\,y^{4} +
0.8066934397\,x\,y - 0.03333333332\,x^{2}\,y  \mbox{} +
0.03999999998\,x^{3}\,y^{3} + 0.5466666665\,x^{2}\,y^{2 } +
0.1133333333\,x^{3}\,y^{2}  \mbox{} + 0.04999999998\,x\,y^{4} -
0.5333333332\,y^{3} + 0.006666666664\,x^{6}  \mbox{} +
0.006666666664\,x^{4}\,y^{2} + 0.006666666664\,y^{6} -
0.6700000000\,x^{2}  \mbox{} - 0.1766666666\,x^{4} +
0.3599999999\,y^{4} - 0.4699999998\,x^{3} - 0.006666666664\,x^{5}
\mbox{} - 0.6066666665\,y^{2} - 0.03999999998\,y^{5}.$ }

\parasmall

 \noindent  We get the $\epsilon$-singular locus $\sing({\cal C})={\cal S}_{1}\cup
{\cal S}_{2}\cup {\cal S}_{3}$ where
\[ \hspace*{-3 mm}
\begin{array}{l}
\begin{array}{ll}
{\cal S}_1=\{\hspace*{-0.3cm}& Q_{1}=(-1.994232333, 1.005043048),\\
          \,& Q_{2}=(-2.000005299+0.005645280797i, -1.000026945-0.0002822677587i),\\
           \,&Q_{3}=(-2.000014217+0.004619269427i,
1.000004775-0.003559494332i),\\
           \,& Q_{4}=(-2.003547061, -1.006293429),
\, Q_{5}=(-2.005740475,
0.9948974977),\\
 \,&Q_{6}=(-1.996418580, -0.9936748962),\\
\, &Q_{7}=
(-2.000014217-.004619269427i, 1.000004775+0.003559494332i), \\
\, &Q_{8}= (-2.000005299-0.005645280797i,
-1.000026945+0.002822677587i), \\
\, &Q_{9}= (1.000036272+0.008596901071i, 2.000017052-0.003059926359i), \\
\, &Q_{10}= (5.999999669, -2.999998564), \\
\, &Q_{11}= (1.000036272-0.008596901071i,
2.000017052-0.003059926359i), \\
\, &Q_{12}=(0.9978910941, 1.994329680),
 \, Q_{13}=(1.002094534, 2.005650021)\},
\end{array}
\\
\begin{array}{ll}
{\cal S}_2=\{\hspace*{-0.3cm}& Q_{14}=(-2.000000001, 1.),
 \, Q_{15}=(-2., -1.000000005),
 \, Q_{16}=(1., 2.)
 \},
\end{array}
\\
\begin{array}{ll}
{\cal S}_3=\emptyset.
\end{array}
\end{array}
\]
The singular cluster decomposition is  (see Fig.
\ref{clusters6}, Left):

 \parasmall

\parbox{15.5cm}{\small  $\cluster_2(Q_{10}) = \{Q_{10}\},\,
\cluster_3(Q_{14})=\{Q_1, Q_3, Q_5, Q_7,
Q_{14}\},\,\\\cluster_3(Q_{15}) =\{Q_2, Q_4, Q_6, Q_8,
Q_{15}\},\,\,\,\mbox{and}\,\,\,
\cluster_3(Q_{16})= \{{Q_9, Q_{11}, Q_{12}, Q_{13}, Q_{16}}\}.$}

\parasmall

\noindent We observe that ${\cal C}$ is $\epsilon$-rational.
 In Step $4$  we obtain  three
$\epsilon$-ramification points:  $P_{1}=(-1.330235522,
\,0.9268173641),  P_{2}=(-1.979908167, \,0.02661222172),$  and
$P_{3}=(-2.700785807, \,-0.07757312293)$. We note that these
points are not in the cluster of each other and they are not in
the clusters of  the $\epsilon$-singularities (see Fig.
\ref{clusters6}, Right).
\begin{center}
\begin{figure}[ht]
\centerline{
\psfig{figure=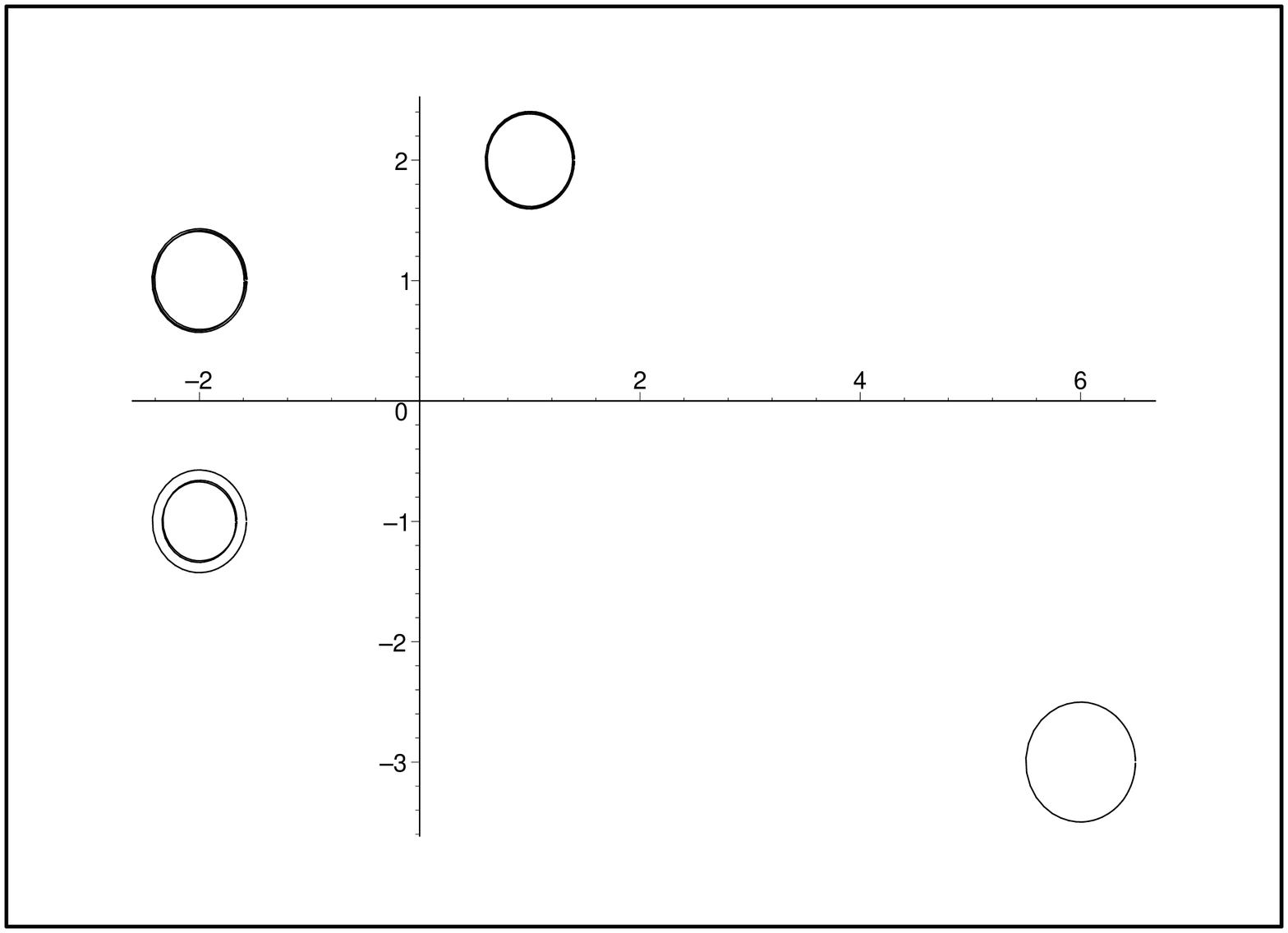,width=5.8cm,height=5.8cm,angle=270} \hspace*{1
cm}
 \psfig{figure=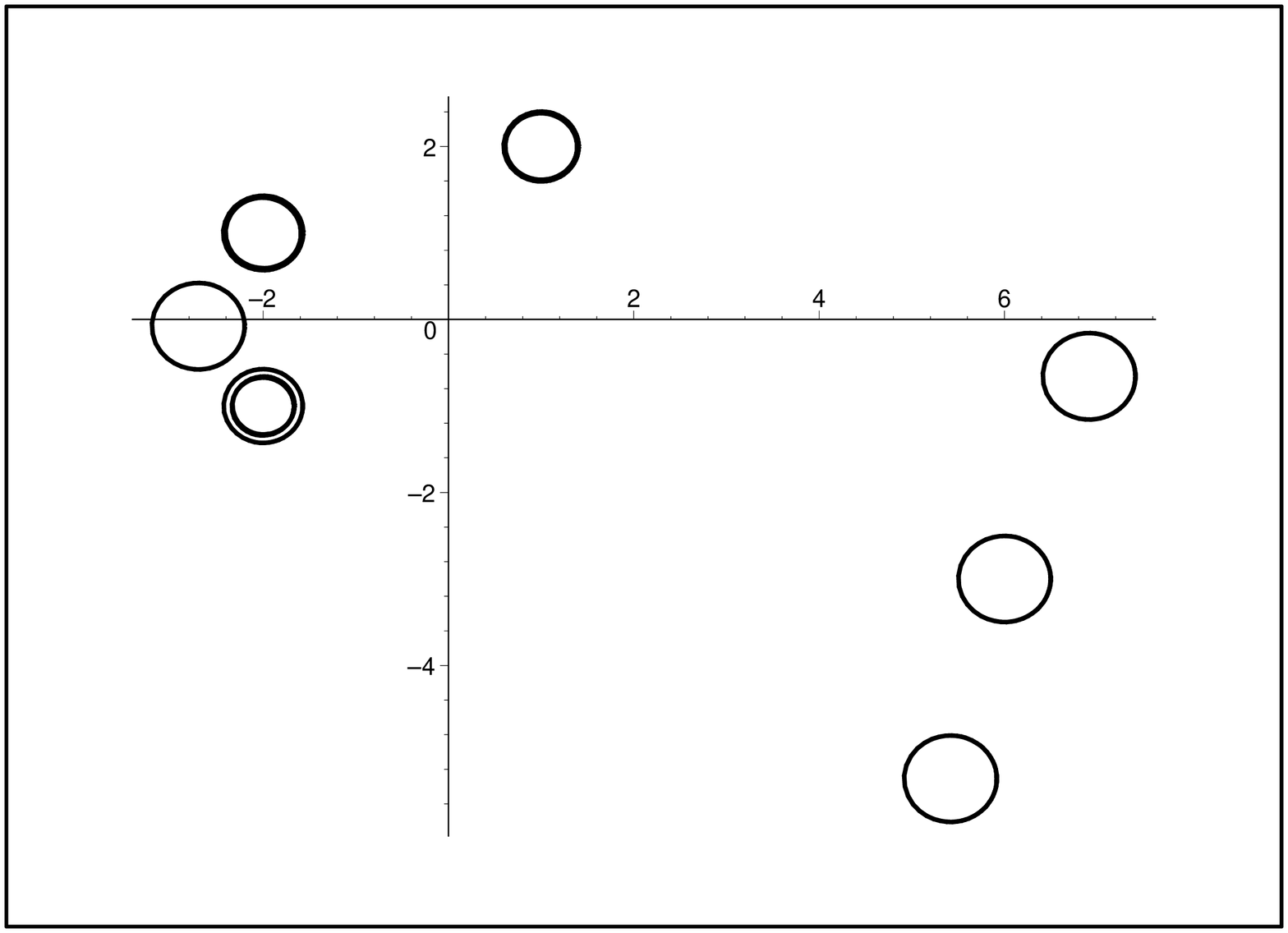,width=5.8cm,height=5.8cm,angle=270
}  } \caption{\small{{\sf Left:} Cluster Decomposition of the
Singular Locus.  {\sf Right:}  Cluster Decomposition of the
Singular Locus with two $\epsilon$-ramification points.}}
\label{clusters6}
\end{figure}
\end{center}
The algorithm  outputs the parametrization
$\overline{\mathcal
P}(t)=(\frac{\overline{p}_1(t)}{\overline{q}(t)},\frac{\overline{p}_2(t)}{\overline{q}(t)})$
 where (see Fig.
\ref{curves6} to compare the input and the output curves):
\begin{center}
\begin{figure}[ht]
\centerline{ \psfig{figure=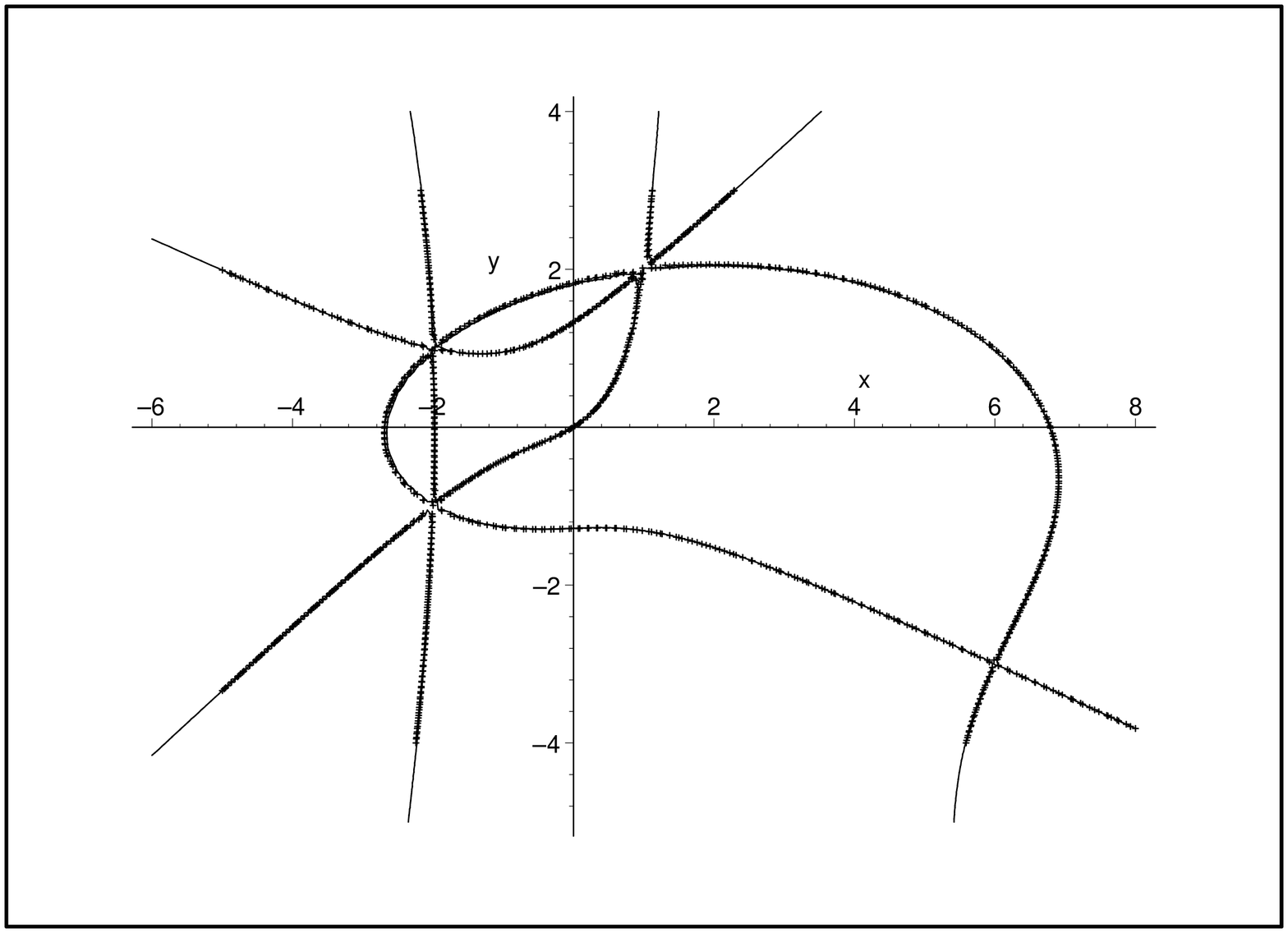,width=6.4cm,height=5.6cm,angle=270}
\psfig{figure=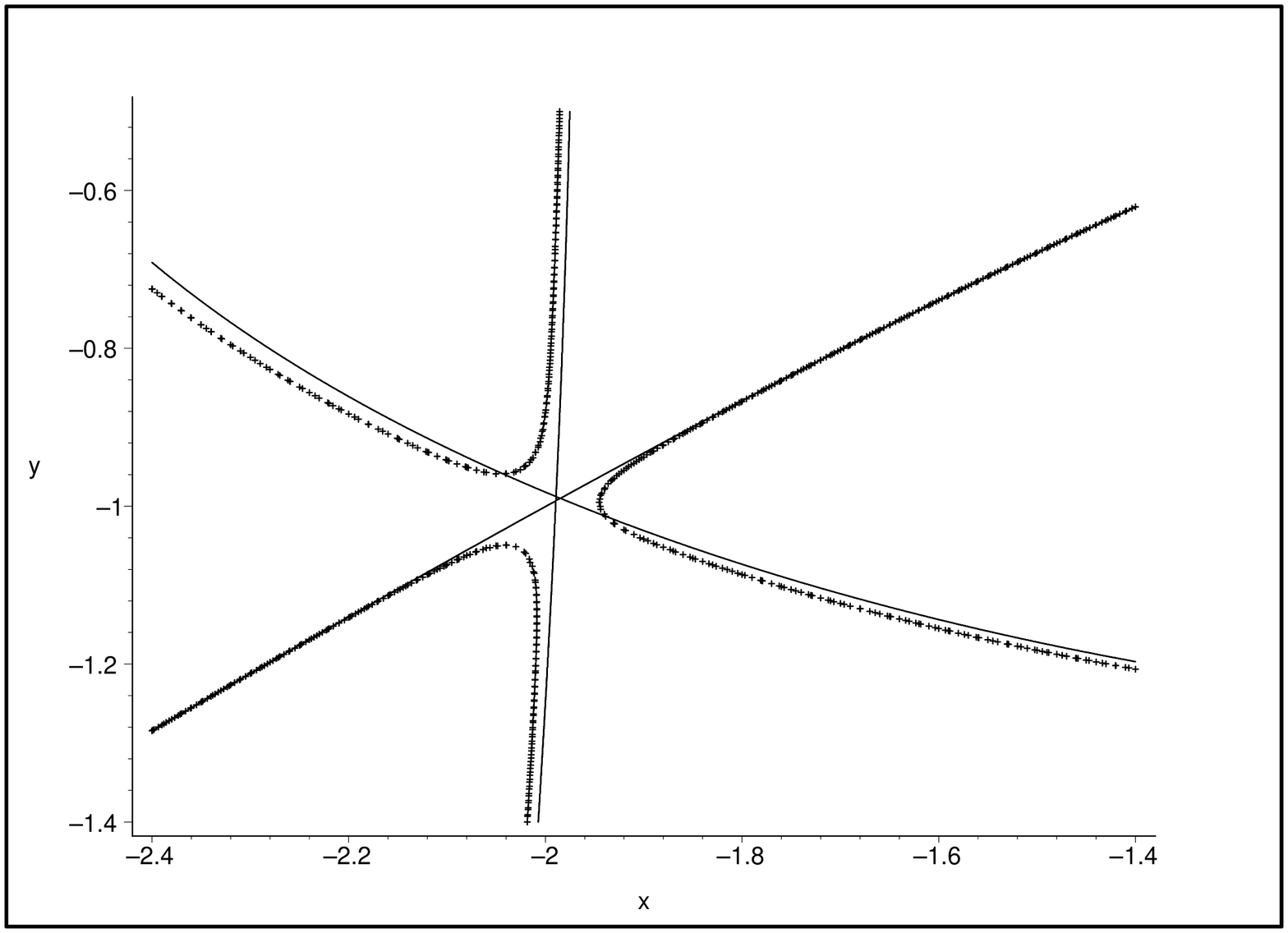,width=6.4cm,height=5.6cm,angle=270}
  } \caption{{\sf Left:} Input (in dots) and output curve in Example \ref{example-degree-6}.
{\sf Right:}  A zoom at $(-2,-1)$
 }
\label{curves6}
\end{figure}
\end{center}

\noindent\parbox{14cm}{\small
 $\overline{p}_1(t)=
 - 0.2992985374\,10^{-13}( - 0.4498780650\,10^{
665}\,t^{2} - 0.1104625259\,10^{666}\,t  \mbox{} +
0.3823432112\,10^{663}\,t^{5} - 0.4349945664\,10^{664} \,t^{3} -
0.8977530140\,10^{665}  \mbox{} + 0.2445532363\,10^{662}\,t^{6} +
0.1487140379\,10^{664} \,t^{4})$}

\noindent
\parbox{14cm}{\small
 $\overline{p}_2(t)= -
0.5410017657\,10^{-14}(0.4254697372\,10^{662 }\,t^{6} -
0.4450231957\,10^{665}\,t^{3}  \mbox{} -
0.8325944623\,10^{665}\,t^{2} + 0.3137371087\,10^{665} \,t -
0.1036628422\,10^{663}\,t^{5} \mbox{} -
0.7148994612\,10^{664}\,t^{4} + 0.1592853294\,10^{666}) $}

\noindent\parbox{14cm}{\small
$\overline{q}(t)=0.3388387927\,10^{649}\,t^{5} \mbox{} +
0.1633128101\,10^{651}\,t^{2} + 0.3082510569\,10^{648} \,t^{6} +
0.4492917291\,10^{651}\,t  \mbox{} + 0.1270749205\,10^{650}\,t^{4} +
0.3531547270\,10^{650} \,t^{3} + 0.4139801407\,10^{651}.$}

\parasmall

 We note that the algorithm did no require
perturbing $\overline{H}^{**}_{d-2}$.\hfill \qed

 \end{example}

 \parasmall

\begin{example}\label{example-degree-7}

Let us consider $\epsilon= 0.001$ and the curve ${\cal C}$ of proper
degree $7$ defined by the polynomial (see Fig.\ref{curves7}):
 \vspace{1 mm}

\parbox{15cm}{\small
 $f(x,y)= 0.005242164122x+0.0000006109092905+0.4234041949y-0.05219720755xy^4-0.1626221914x^3y^2-
 0.006150324474x^5y^2-0.009115378696y^3+0.01468749412x^3-0.1726592957y^4-0.005178781717x^5+
 0.0006102983812y^6+0.7056394692y^2+0.007271579029x^4-0.009049345878x^3y^4+0.02810421594y^5+
 0.01517536020x^4y^3-0.03335531981xy^3+0.07030423460x^3y^3+0.9999999999xy^2+0.02396026447x^5y
 +0.06359239287x^2y^3+0.0006102983812xy^5+0.06037915453x^4y-0.05961614786x^4y^2+0.1735938027x^2y^2+
 0.009673386920xy-0.1183998159x^3y-0.3997312415x^2y-
 0.01504641433x^2y^4-0.0002034043985x^7+0.0007152781730x^6y+0.009647777670xy^6-0.01996027092x^2-
 0.001858780227x^6-0.008636725103x^2y^5-0.002554427076y^7.$ }

\parasmall

\noindent The  $\epsilon$-singular locus is $\sing({\cal C})={\cal S}_{1}\cup
{\cal S}_{2}\cup {\cal S}_{3}$ where
\[ \hspace*{-3 mm}
\begin{array}{l}
\begin{array}{ll}
{\cal S}_1=\{\hspace*{-0.3cm}& Q_{1}=(4.998181206+0.0004639222080i, 6.997094116-0.0003357295061i),\\
          \,& Q_{2}=(4.998181206-0.0004639222080i, 6.997094116+0.0003357295061i),\\
           \,&Q_{3}=(5.001816967+0.0004635470406i,
7.002902744-0.0003352187676i),\\
           \,& Q_{4}=(0.9999998537, -3.000000118),\\
\,&Q_{5}=(5.001816967-0.0004635470406i,
7.002902744+0.0003352187676i),
\\
 \,&Q_{6}=(-2.000211362+0.00008683312445i,
-0.001218550314-0.9986341029i),\\
\, &Q_{7}= (-2.000211362-0.00008683312445i,
-0.001218550314+0.9986341029i), \\
\, &Q_{8}= (1.998594026+0.0001453051485i,
-0.0005005279912-0.9994646423i), \\
\, &Q_{9}= (-1.998594026-0.0001453051485i,
-0.0005005279912+0.9994646423i),\\
\,\,\,\,\, &Q_{10}= (1.000001333,
-0.00000005539458512), \\
\, &Q_{11}= (-2.001405450-0.0001446569643i,
0.0005000190394-1.000535463i), \\
\, &Q_{12}= (-2.001405450+0.0001446569643i,
0.0005000190394+1.000535463i), \\
\, &Q_{13}=(-1.999787395-0.00008416274464i,
0.001216340837-1.001365249i), \\
\, &Q_{14}=(-1.999787395+0.00008416274464i, 0.001216340837+1.001365249i), \\
\, &Q_{15}=(4.997608917-0.001980994691i,
1.999734804-0.002346469999i), \\
\, &Q_{16}=(4.997608917+0.001980994691i,
1.999734804+0.002346469999i),\\
\, &Q_{17}=(-3.999997183, 1.999998082),\\
\, &Q_{18}=(5.002393988-0.001973712849i, 2.000267815-0.002341789270i),\\
 \, &Q_{19}=(5.002393988+0.001973712849i, 2.000267815+0.002341789270i)\},
\end{array}
\\
\begin{array}{ll}
{\cal S}_2=\{\hspace*{-0.3cm}&  Q_{20}=(-2.000000398-0.0000003109941563i, \\ & \hspace*{1.2 cm} -0.0000005243819124-0.9999997083i), \\
 \, &Q_{21}=(-2.000000398+0.0000003109941563i, \\ & \hspace*{1.2 cm}-0.0000005243819124+0.9999997083i),\\
 \, &Q_{22}=(5.000000495, 2.000000179),
 \, Q_{23}=(4.999999480, 6.999999337)
 \},
\end{array}
\\
\begin{array}{ll}
{\cal S}_3=\emptyset.
\end{array}
\end{array}
\]

\noindent The cluster decomposition of the singular locus is
(see Fig. \ref{clusters6}, Left):

 \parasmall

\parbox{15.5cm}{\small  $\cluster_2(Q_{4}) = \{Q_{4}\},\,
\cluster_2(Q_{10}) = \{Q_{10}\},\,\cluster_2(Q_{17}) =
\{Q_{17}\},\,\\
 \cluster_3(Q_{23})=\{Q_1, Q_2, Q_3, Q_5,
Q_{23}\},\,\,\,\cluster_3(Q_{20}) =\{Q_6, Q_8, Q_{11}, Q_{18},
Q_{20}\}, \\
\cluster_3(Q_{21})= \{{Q_7, Q_{9}, Q_{12}, Q_{14},
Q_{21}}\},\,\,\,\mbox{and}\,\,\,\cluster_3(Q_{22})= \{{Q_{15},
Q_{16}, Q_{18}, Q_{19}, Q_{22}}\}.$}

\parasmall

\noindent We observe that ${\cal C}$ is $\epsilon$-rational.
 Now,  we obtain  four
$\epsilon$-ramification points:  $P_{1}=( -2.972405737,
-7.933174980),$  $ P_{2}=( 23.79950366, 17.84891277),$
$P_{3}=(-10.06218879, 1.300686562) $ and $P_{4}=(24.47385001,
17.37936091)$. We note that these  points are not in the cluster
of each other and they are not in the clusters of  the
$\epsilon$-singularities.
 Finally,  the algorithm   outputs a parametrization (for space limitation we do not include it here).
 In Fig. \ref{curves7} we plot the input and the input curve.
\begin{center}
\begin{figure}[ht]
\centerline{ \psfig{figure=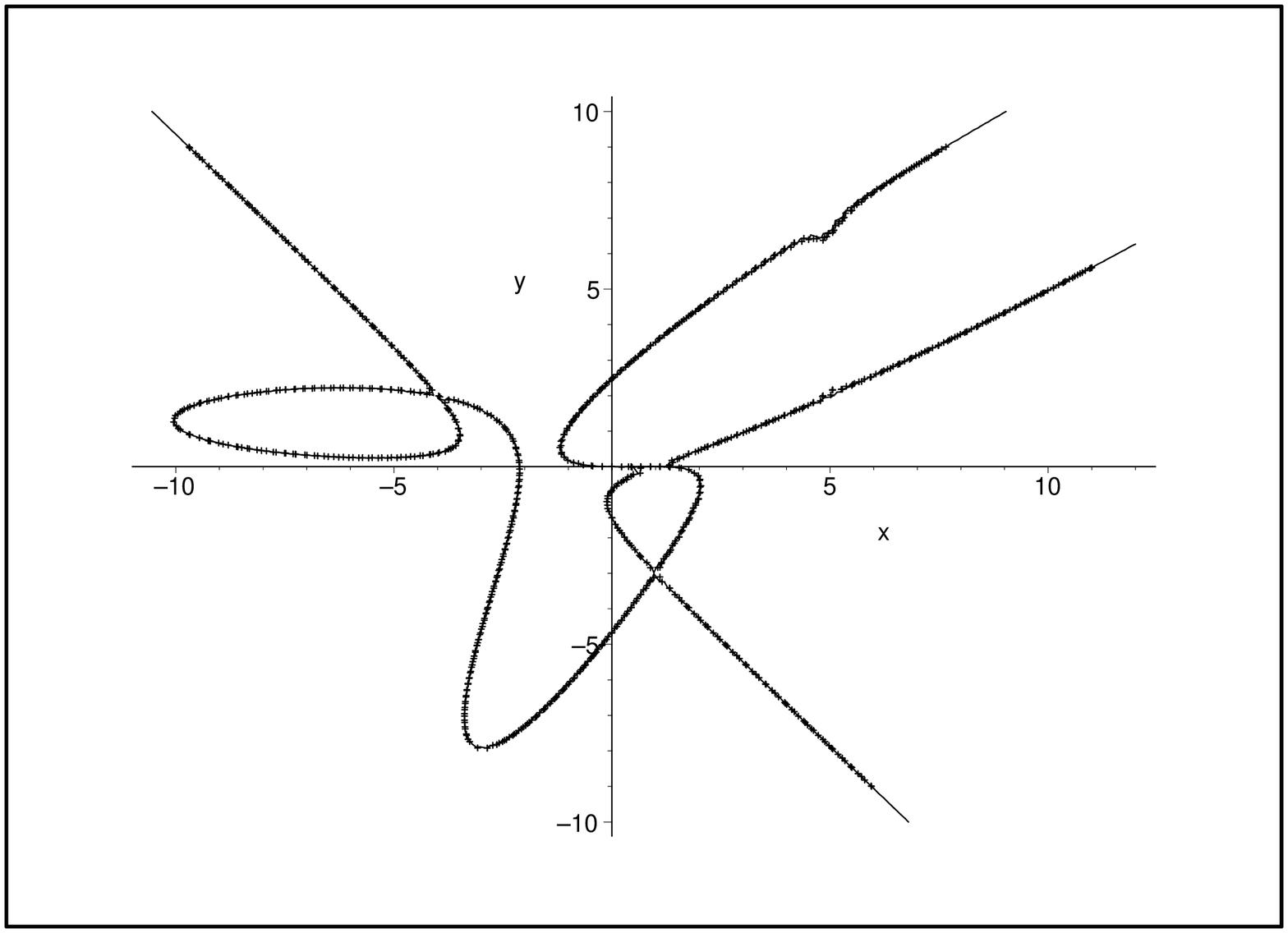,width=6.4cm,height=5.6cm,angle=270}
\psfig{figure=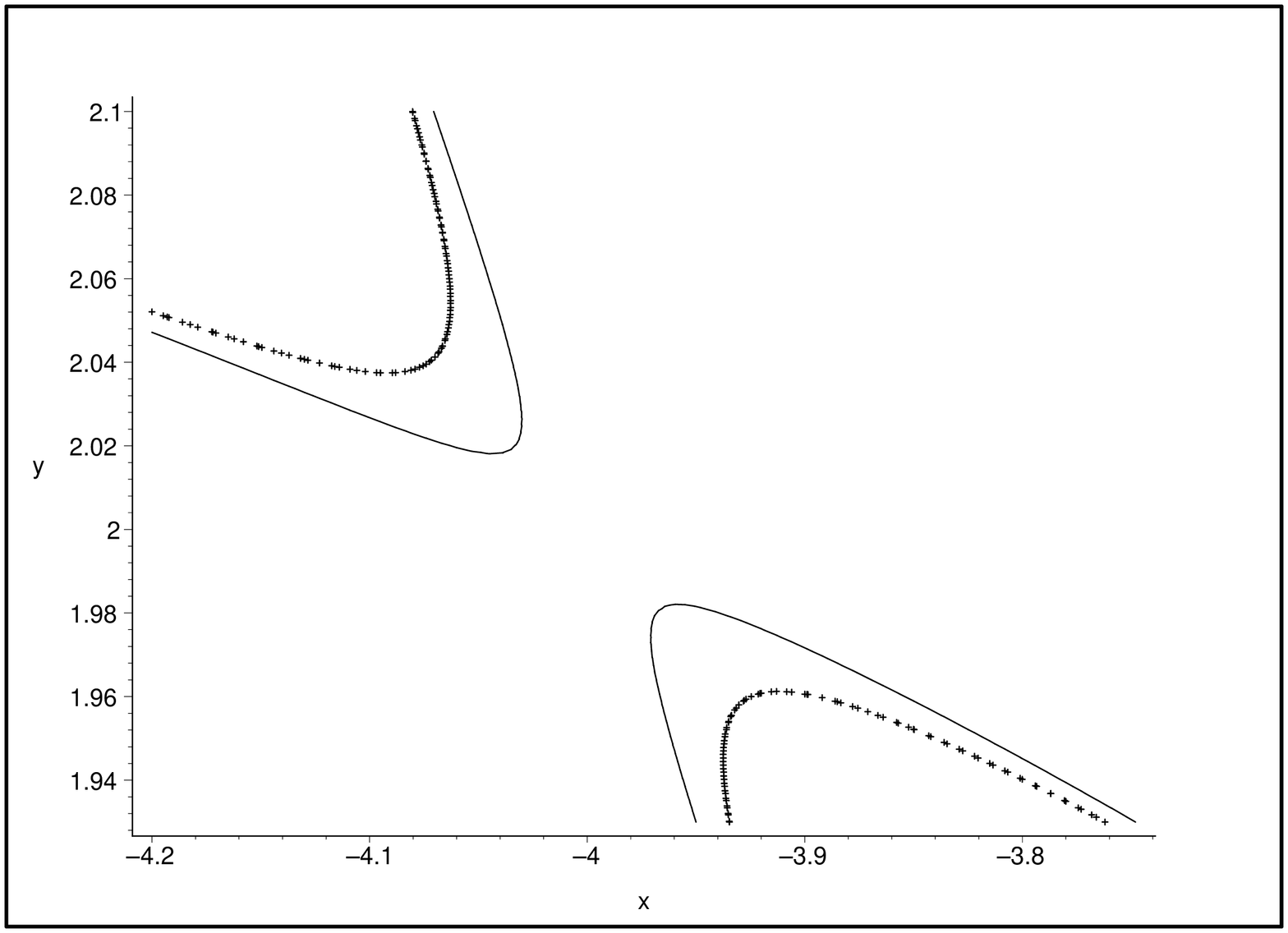,width=6.4cm,height=5.6cm,angle=270}
  } \caption{{\sf Left:} Input (in dots) and output curve in Example \ref{example-degree-7}.
  {\sf Right:}  A zoom at $(-4,2)$
 }
\label{curves7}
\end{figure}
\end{center}
\vspace*{-1.5cm}
\hfill \qed

\end{example}

\subsection{Empirical Analysis of Error}

\begin{center}
\begin{figure}[ht]
\centerline{
\psfig{figure=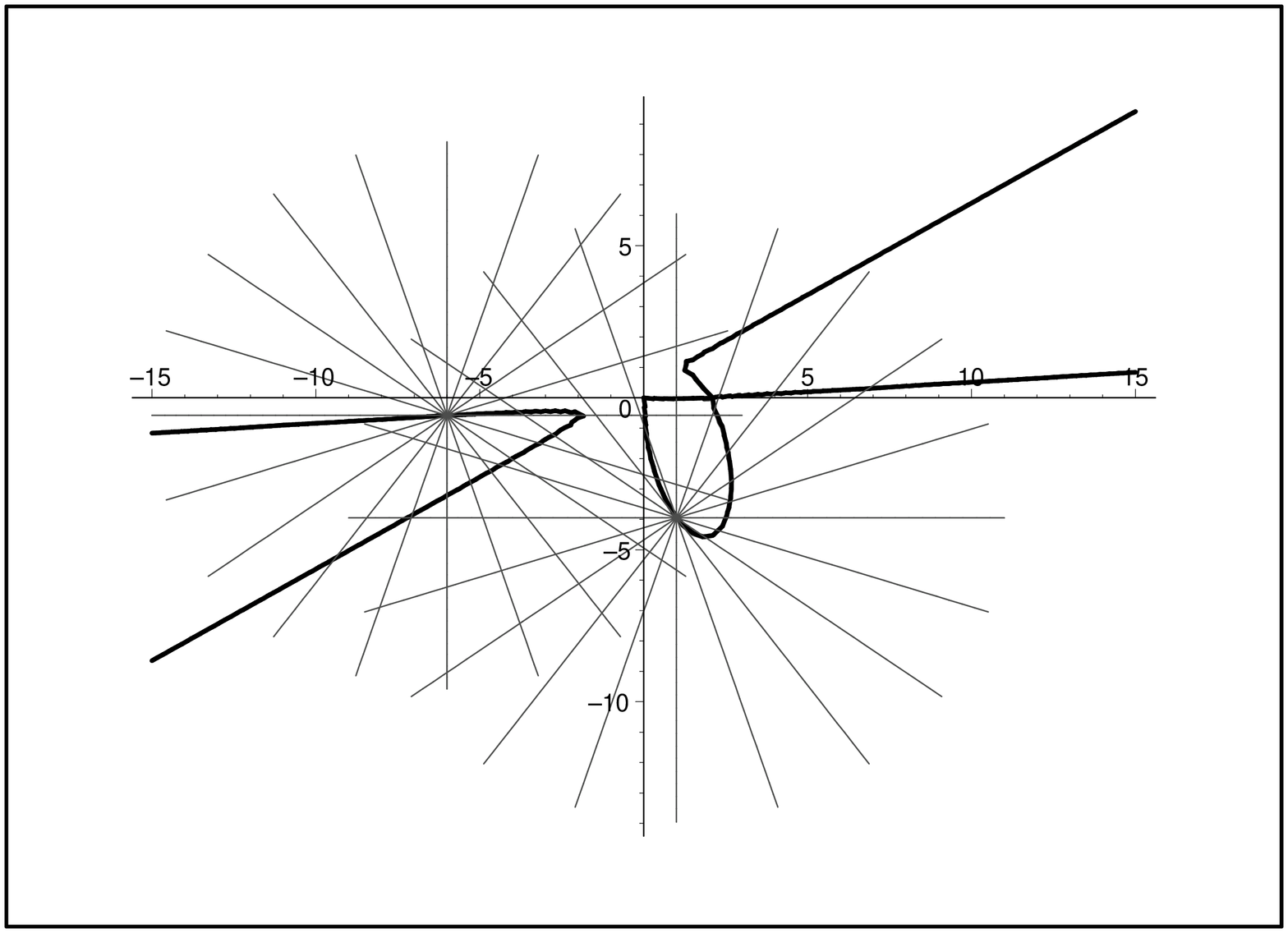,width=6.5cm,height=5.5cm,angle=270}
\psfig{figure=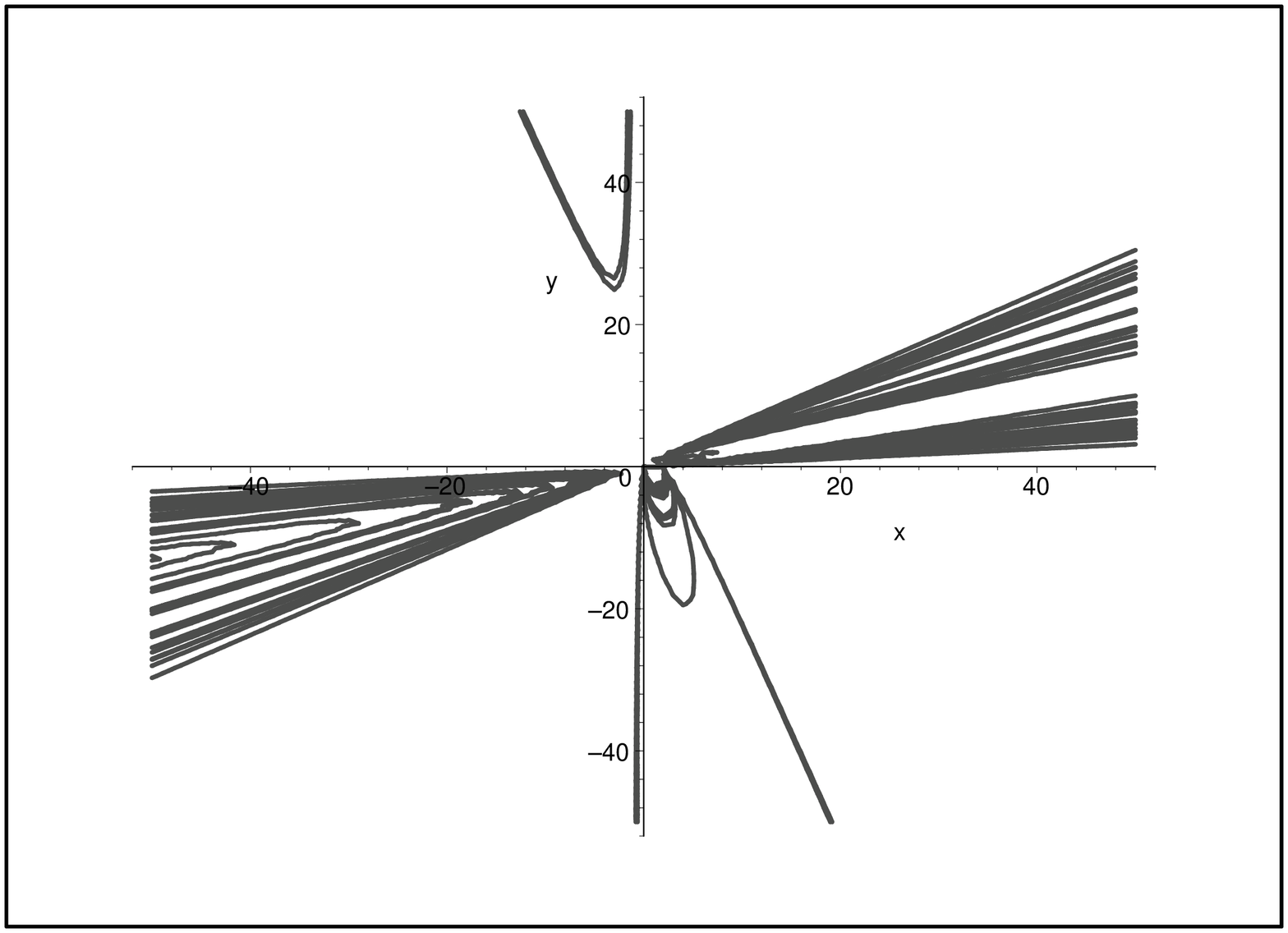,width=6.5cm,height=5.5cm,angle=270}
 }
\caption{ {\sf Left:} Illustration of the method. {\sf Right:} The
28 $\epsilon$--rational curves ${\cal C}_{ij}$ randomly generated
. } \label{fig_fam}
\end{figure}
\end{center}

The aim of this section is to analyze empirically the performance
of the algorithm proposed. A good performance would mean to obtain
an output curve $\overline{\cal C}$ close to the input curve $\cal
C$ and by close we mean that  $\cal C$ is contained in the offset
region of $\overline{\cal C}$ at a small distance and viceversa
(see \cite{PSS}). To estimate the distance, between the curves
$\cal C$ and $\overline{\cal C}$, we designed the next method (see
Fig. \ref{fig_fam}, Left):
\begin{enumerate}
\item We randomly generate a set $\cal E$ of (affine)
$\epsilon$--points on the input curve $\cal C$ as follows. Fix
real numbers $a,b$, and take $n$ random integer values
$\alpha_i\in [a,b]$, $i=1,\ldots ,n$. Let ${\cal E}_1$ be the set
of intersection points of the curve $\cal C$ with the  lines
$x=\alpha_i$. We also take $n$ random integer values $\beta_i\in
[a,b]$ and we intersect the curve $\cal C$ with the lines
$y=\alpha_i$ to obtain the set of points ${\cal E}_2$. We set
${\cal E}={\cal E}_1\cup{\cal E}_2$. \item Let $r$ be a positive
integer and $\Theta_r=\{\frac{k \pi}{r}\mid k=1,\ldots ,r\}$. For
each $\epsilon$--point $P$ in $\cal E$ and each
$\theta\in\Theta_r$ let $L_{P,\theta}$ be the line through $P$ in
the direction of $(cos(\theta),sin(\theta))$.

\item For each $P\in{\cal E}$ we compute $d_{P,\overline{\cal
C}}=\min \{\|P-Q\|_2 \mid Q\in \cup_{\theta\in\Theta_r}
(L_{P,\theta}\cap \overline{\cal C})\}.$

\item Let ${\cal D}=\{d_{P,\overline{\cal C}}\mid P\in{\cal E}\}$.
We compute the mean value $\mu$ of the elements of ${\cal D}$ as
well as the statistical standard error $\rho$. We can say that the
estimated distance is, in average, in the interval
$[\mu-1.96\,\rho,\mu+1.96\,\rho]$.
\end{enumerate}

\noindent Given the curve ${\cal C}$ of degree $d=4$ defined by the
irreducible polynomial
\begin{equation}\label{eq_error_curve}
\begin{array}{l}f(x,y)=1.000065 y^2+1.00000028 y^3+y^4+1.000065 xy
-11.49999972xy^2  \\ +xy^3  +0.760065 x^2+5.74000028 x^2y
 +3.69x^2y^2-0.75999972x^3-3.12x^3y
\\ +0.19x^4 +0.01x+0.01y.
\end{array}
\end{equation}
For $\epsilon=0.01$ the algorithm concluded that ${\cal C}$ is
$\epsilon$-rational and returned a rational parametrization ${\cal
P} (t)=(\overline{p}_1(t),\overline{p}_2(t))$ which corresponds to
the rational curve $\overline{\cal C}$ with implicit equation
\begin{equation*}
\begin{array}{l}\overline{f}(x,y)=0.01642553x^4+0.06494377x^2+0.08804654y^2-0.06552116x^3\\
+0.08169391y^3+0.091025077xy+0.49976135x^2y-0.99999999xy^2+0.08645018y^4\\
+0.31900118x^2y^2-0.26972458x^3y+0.08645019xy^3-0.00001398+0.00078781x\\
+0.0007408y.
\end{array}
\end{equation*}
With the above notation, we take $[a,b]=[-100,100]$, $n=15$. The
number of points used to compute the experimental distance is $|
{\cal E} |=120$. Set the number of lines going through each point
equal $r=10$. The mean value of the positive integers in $\cal D$
is $\mu=0.007541$ and the statistical standard error
$\rho=0.000855$. We can conclude that the estimated distance is,
in average, in the interval
$[\mu-1.96\,\rho,\mu+1.96\,\rho]=[0.005866, 0.009217]$.

For the examples treated in Section \ref{sec-examples},  the
corresponding figures show  that input and output curves  look
very close. We show below that the estimated distance between
$\overline{\cal C}$ and $\cal C$ is small compared to the value of
$\epsilon$ used in each one of the examples. For this purpose, we
execute the previously described method to estimate the distance
between $\overline{\cal C}$ and $\cal C$ taking $r=10$, and
$[a,b]=[-100,100]$. Let $\mid {\cal E}\mid$ be the number of
points taken to compute the experimental distance. The results
obtained in the examples are:
\begin{enumerate}
\item Let $\cal C$ and $\overline{\cal C}$ be the curves in
Example \ref{example-degree-5}, $\epsilon=0,01$, $\mid {\cal
E}\mid=150$. Then, $\mu= 0.070953, \rho=0.012655$. So, the
estimated distance is, in average, in the interval
$[0.046149,0.095757].$ \item Let $\cal C$ and $\overline{\cal C}$
be the curves in Example \ref{example-degree-6}, $\epsilon=0,004$,
$\mid {\cal E}\mid=180$. Then, $\mu=0.002425, \rho=0.000100$. So,
the estimated distance is, in average, in the interval $[0.002228,
0.002621].$ \item Let $\cal C$ and $\overline{\cal C}$ be the
curves in Example \ref{example-degree-7}, $\epsilon=0,001$, $\mid
{\cal E}\mid=210$. Then, $ \mu=0.006560, \rho=0.000609$. So, the
estimated distance is, in average, in the interval $[0.005366,
0.007754].$
\end{enumerate}
Until now we have empirically measured the distance in the
examples used in paper. We describe next a different experiment.
We randomly generate a set of curves and for each curve we
estimate its distance to the output curve given by our algorithm.
Our experiments are satisfactory and allow us to think about a
theoretical treatment of this fact as a future project. We explain
next how the family of curves was constructed. We fix three points
$P_1=(2:0:1)$,$P_2=(0:0:1)$ and $P_3=(1:1:1)$ in
$\mathbb{P}^2({\Bbb C})$. We consider the linear system of curves
of degree 4 defined by the divisor $2P_1+2P_2+2P_3$. Its defining
polynomial is
\begin{equation}\label{G}
\begin{array}{l}G(x,y,z)=u_{2}y^2z^2+u_{3}y^3z+u_{4}y^4+u_{5}xyz^2+(-2u_{2}-3u_{3}\\
-4u_{4}-1/2u_{5}-2u_{6})xy^2z+u_{6}xy^3+u_{1}x^2z^2+(-3/2u_{5}+2u_{3}\\
+4u_{4}+2u_{6}-u_{1})x^2yz+(u_{2}+u_{3}+1/2u_{5}+1/4u_{1}+u_{4})x^2y^2\\
-u_{1}x^3z+(1/2u_{5}-u_{3}-2u_{4}-u_{6}+1/2u_{1})x^3y+1/4u_{1}x^4
\end{array}
\end{equation}

For $j=1,\ldots ,6$ and $i=1,\ldots ,10$ let $r_{ij}$ be a random
integer number in the interval $[0,100]$. We obtain 60 different
polynomials $G_{ij}(x,y,z)$, $j=1,\ldots, 6$, $i=1,\ldots ,10$
 setting
\[u_k=\left\{\begin{array}{lcc}(\frac{r_{ij}}{100})^i&\mbox{ if }&k=j\\ 1&\mbox{ if }&k\neq j\end{array}\right. \,\,\, k=1,\ldots ,6\]
in equation ~\eqref{G}.
Given $i\in\{1,\ldots ,6\}$ and $j\in \{1,\ldots ,10\}$ we obtain a random perturbation $g_{ij}(x,y)\in \mathbb{R}[x,y]$ of $G_{ij}(x,y,1)$ as follows
\[g_{ij}(x,y)= G_{ij}(x,y,1)+ \epsilon \frac{r_1}{100} (x+y)+\epsilon^2 \frac{r_2}{100}(x^2+xy+y^2)+
\epsilon^3 \frac{r_3}{100}(x^3+x^2y+xy^2+y^3)\] where
$r_1,r_2,r_3$ are integer numbers taken randomly in the interval
$[0,100]$ and $\epsilon=0.01$. The polynomials $g_{ij}(x,y)$,
$j=1,\ldots, 6$, $i=1,\ldots ,10$ have proper degree 4 and define
60 curves ${\cal C}_{ij}$ verifying $(1:0:0),(0:1:0)\notin {\cal
C}_{ij}^h$  and without (exact) singularities at infinity. Using
our algorithm  we  conclude that 28 of the 60 curves are
$\epsilon$--rational. We show those curves in Fig. \ref{fig_fam},
Right. The implicit equation of the curve ${\cal C}_{11}$ is
equation \eqref{eq_error_curve}.

With the notation used in the method described above to compute
the experimental distance, we take $[a,b]=[-100,100]$ and $n=15$
so that the number of points used to compute the distance is
expected to be $\mid {\cal E}\mid=120$. Set the number of lines
going through each point equal $r=10$. We compute $\mu$,$\rho$ and
$I_{\mu,\rho}=[\mu-1.96\,\rho,\mu+1.96\,\rho]$ for each one of the
28 $\epsilon$--rational curves to obtain the next table.

\small
\[ \hspace*{-2 mm}
\begin{array}{cc}
\mbox{\small\begin{tabular}{|c||c||c|}
\hline
$\mu$ & $\rho$ & $I_{\mu,\rho}$\\
\hline
$0.007541$ & $0.000855 $& $[0.005866, 0.009217]$
\\
\hline $0.006977$& $0.001184$&$[0.004656, 0.009299]$
\\
\hline
0.006977& 0.001184&[0.004656, 0.009299]\\
\hline
0.003577& 0.000503&[0.002592, 0.004563]\\
\hline
0.004011& 0.000553&[0.002928, 0.005094]\\
\hline
0.006007& 0.000808&[0.004423, 0.007590]\\
\hline
0.004239& 0.000844&[0.002584, 0.005894]\\
\hline
0.005758& 0.000585&[0.004610, 0.006905]\\
\hline
0.002882& 0.000224&[0.002442, 0.003322]\\
\hline
0.005477& 0.000756&[0.003996, 0.006958]\\
\hline
0.003123& 0.000437&[0.002266, 0.003979]\\
\hline
0.004752& 0.000359&[0.004049, 0.005455]\\
\hline
0.001453& 0.000148&[0.001163, 0.001744]\\
\hline
0.004956& 0.007123&[0.035599, 0.063522]\\
\hline
0.001049& 0.000113&[0.000827, 0.001272]\\
\hline
\end{tabular}} &
\mbox{\small\begin{tabular}{|c||c||c|}
\hline
$\mu$ & $\rho$ & $I_{\mu,\rho}$\\
\hline
0.006807& 0.000385&[0.006051, 0.007563] \\
\hline
0.100902& 0.013253&[0.074926, 0.126879]\\
\hline
0.003049& 0.000254&[0.002551, 0.000355]\\
\hline
0.003924& 0.000212&[0.003508, 0.004339]\\
\hline
0.003995& 0.000549&[0.002919, 0.005072]\\
\hline
0.008330& 0.000806&[0.006749, 0.009911]\\
\hline
0.005638& 0.000536&[0.004586, 0.006690]\\
\hline
0.003020& 0.000316&[0.002399, 0.003639]\\
\hline
0.000854& 0.000091&[0.000677, 0.001032]\\
\hline
0.004077& 0.000274&[0.003540, 0.004614]\\
\hline
0.035130& 0.005220&[0.024898, 0.045361]\\
\hline
0.006209& 0.000619&[0.004996, 0.007423]\\
\hline
0.013406& 0.001179&[0.011094, 0.015718]\\
\hline
 0.009037&0.0006872&[0.007691, 0.010385]\\
\hline
\end{tabular}}
\end{array}
\]

\end{document}